\documentclass[10pt]{article}[10pt]

\usepackage[centertags]{amsmath}
\usepackage{amssymb,amsfonts,theorem,stmaryrd}
\usepackage[all]{xy}

\allowdisplaybreaks

\newtheorem{thm}{Theorem}[section]
\newtheorem{prop}[thm]{Proposition}

\newtheorem{lemma}[thm]{Lemma}
\newtheorem{cor}[thm]{Corollary}
\theorembodyfont{\rmfamily}

\newenvironment{pf}{\textit{Proof.}}{\hfill$\boxbox$}

\newcommand{\F}{\mathcal{F}}
\newcommand{\gM}{g^{TM}}
\newcommand{\hF}{h^{\mathcal{F}}}
\newcommand{\nF}{\nabla ^\mathcal{F}}
\newcommand{\dF}{d(\mathcal{F})}

\newcommand{\lF}{\Delta(\mathcal{F})}
\newcommand{\Tr}{\operatorname{Tr}}
\newcommand{\Pf}{\operatorname{Pf}}
\newcommand{\Id}{\operatorname{Id}}
\newcommand{\Str}{\operatorname{Str}}
\newcommand{\M}{\operatorname{\mathsf{M}}}
\newcommand{\E}{\operatorname{\mathbb{E}}}
\newcommand{\LIM}{\operatorname*{\mathsf{LIM}}}
\newcommand{\Hom}{\operatorname{Hom}}
\newcommand{\End}{\operatorname{End}}
\newcommand{\Gr}{\operatorname{Gr}}
\newcommand{\Cas}{\operatorname{Cas}}
\newcommand{\T}{\operatorname{T}}
\newcommand{\N}{\operatorname{N}}

\newcommand{\Irr}{\operatorname{Irr}}

\begin{document}
%%%%%%%%%
% title %
%%%%%%%%%
\title{The variation formulas for the equivariant Ray-Singer metric}
\author
{Hartmut Wei{\ss} \\ LMU M\"unchen}

\maketitle

%%%%%%%%%%%%
% abstract %
%%%%%%%%%%%%
\begin{abstract}
We give a new and detailed proof of the variation formulas for the equivariant
Ray-Singer metric, which are originally due to J.M.~Bismut and W.~Zhang. 
\end{abstract}

%%%%%%%%%%%%%%%%%%%%%
% table of contents %
%%%%%%%%%%%%%%%%%%%%%
\tableofcontents

%%%%%%%%%%%%%%%%
% Introduction %
%%%%%%%%%%%%%%%%
\section{Introduction}
Let $M$ be a closed $n$-dimensional smooth manifold and $(\mathcal{F},\nabla^{\mathcal{F}})$ a flat complex vector bundle over $M$. Let further $g^{TM}$ be a Riemannian metric on $M$ and $\hF$ a hermitian metric on $\F$. We will not assume $\hF$ to be parallel with respect to $\nF$.
Let $G$ be a compact Lie group acting smoothly on $M$ such that the metrics $\gM$, $\hF$ and the flat connection $\nF$ are preserved.

To this data one associates the equivariant Ray-Singer analytic torsion $\tau(M,\mathcal{F};\gM,\hF)$ and the equivariant Ray-Singer metric
\begin{equation*}
\|\cdot\|_{\det (H^{\bullet}(M,\mathcal{F}),G)}=|\cdot|_{\det (H^{\bullet}(M,\mathcal{F}),G)}\cdot\tau\bigl(M,\mathcal{F};\gM,\hF\bigr)
\end{equation*}
on the equivariant determinant of $H^{\bullet}(M,\mathcal{F})$, where $|\cdot|_{\det (H^{\bullet}(M,\mathcal{F}),G)}$ is the equivariant $L^2$-metric on $\det (H^{\bullet}(M,\mathcal{F}),G)$. For details we refer the reader to Section \ref{equiv_anal_torsion} and to the original work of J.M.~Bismut and W.~Zhang as published in \cite{BZ1} and \cite {BZ2}. 

An obvious question is to what extent these quantities depend on the geometric data, i.e.~the metrics $\gM$ and $\hF$. The aim of this article is to give a detailed proof of the following result, which is the differential version of Theorem 2.7 in \cite{BZ2}:

\begin{thm}[J.-M.~Bismut, W.~Zhang]\label{main}
For $\gamma \in G$ one has for the variation of the equivariant Ray-Singer metric:\\
\\
(1) $\varepsilon \mapsto\gM(\varepsilon)$:
\begin{equation*}
\frac{\partial}{\partial
 \varepsilon}\Bigr\vert_{\varepsilon=0}\log\|\cdot\|_{(\det
H^{\bullet}(M,\mathcal{F}),G)}^2(\gamma)=-\int\limits_{M^{\gamma}}\theta(\gamma,\mathcal{F},\hF)\wedge\tilde{e}^{\prime}(TM^{\gamma})
\end{equation*}
(2) $\varepsilon \mapsto\hF(\varepsilon)$:
\begin{equation*}
\frac{\partial}{\partial
 \varepsilon}\Bigr\vert_{\varepsilon=0}\log\|\cdot\|_{(\det
H^{\bullet}(M,\mathcal{F}),G)}^2(\gamma)=\int\limits_{M^{\gamma}}\Tr[\gamma^{\mathcal{F}} V]e(TM^{\gamma},\nabla^{TM^{\gamma}})
\end{equation*}
Here $V=(\hF)^{-1}\dot{h}^{\F}$ and $M^\gamma$ denotes the fixed point set of $\gamma \in G$.
\end{thm}
For the defintion of the Euler form $e(TM^{\gamma},\nabla^{TM^{\gamma}})$, the
transgression form $\tilde{e}^{\prime}(TM^{\gamma})$ and the 1-form
$\theta(\gamma,\mathcal{F},\hF)$ we refer the reader again to Section
\ref{equiv_anal_torsion}. Note that the fixed point set $M^\gamma$ is a
compact submanifold without boundary, cf.~\cite{Kob}.

In \cite{BZ1}, J.M.~Bismut and W.~Zhang prove the non-equivariant version of Theorem \ref{main} using a variant of the Getzler rescaling technique, whereas they do not give details in the equivariant case. Our proof is modelled on the proof of the (local) equivariant index theorem by N.~Berline and M.~Vergne in \cite{BV}, see also \cite{BGV}.

{\bf Acknowledgements.} This article is based on the author's diploma thesis written under the supervision of Prof.~Dr.~Ulrich Bunke. The author would like to thank Ulrich Bunke for his support.

%%%%%%%%%%%%%%%%%%%%%%%%%%%%%%%%%%%%%
% The equivariant Ray-Singer metric %
%%%%%%%%%%%%%%%%%%%%%%%%%%%%%%%%%%%%%
\section{The equivariant Ray-Singer metric}\label{equiv_anal_torsion}

Let $\mathcal{A}^\bullet(M,\mathcal{F})=\Gamma(M,\Lambda^\bullet{T}^\ast{M}\otimes\mathcal{F})$ denote the differential forms on $M$ with values in $\F$. Let further
\begin{equation*}
\dF{}:{}\mathcal{A}^\bullet(M,\mathcal{F})\longrightarrow\mathcal{A}^{\bullet+1}(M,\mathcal{F})
\end{equation*}
denote the exterior differential associated with the flat connection $\nF$. The Hodge Laplacian is given by $\lF=\dF\dF^*+\dF^*\dF$, where $\dF^*$ denotes the formal adjoint of $\dF$. For $t >0$ let $\exp(-t\lF)$ denote the heat operator.

For $\gamma \in G$ and $s \in \mathbb{C}
$ let
\begin{equation*}
\theta(\gamma,s)=\frac{1}{\Gamma(s)}\int_{0}^{\infty}t^{s-1}\Tr\bigl\{(-1)^{\N}\N\gamma\exp(-t\lF)(1-P_0)\bigr\}
dt\,,
\end{equation*}
where
 $\N:\mathcal{A}^{\bullet}(M,\mathcal{F})\rightarrow\mathcal{A}^{\bullet}(M,\mathcal{F})$
 is the number operator, which multiplies a homogeneous form with its degree,
 and $P_0$ is the harmonic projection. Hence $\theta(\gamma,s)$ is
 the Mellin transform of
$f(t)=\Tr\bigl\{(-1)^{\N}\N\gamma\exp(-t\lF)(1-P_0)\bigr\}$, i.e.
\begin{equation*}
\M[f](s)=\frac{1}{\Gamma(s)}\int_{0}^{\infty}t^{s-1}f(t)dt
\end{equation*}
Using the asymptotic expansion of the heat kernel one shows that $s \mapsto \theta(\gamma,s)$ is a meromorphic function on the complex plane which is holomorphic about $s=0$. Therefore the {\em equivariant Ray-Singer analytic torsion}
\begin{equation*}
\tau\bigl(M,\mathcal{F};\gM,\hF\bigr)(\gamma):=\exp\Bigl(-\frac{1}{2}\,\frac{\partial}{\partial
s}\Bigr\vert_{s=0}\theta(\gamma,s)\Bigr)
\end{equation*}
is well defined.

Let $V$ be a complex $G$-representation. Let $\det V$ denote the highest exterior power of $V$. We consider the isotypical decomposition
$$
V = \bigoplus_{W \in \hat{G}}\Hom_G(W,V) \otimes_{\mathbb{C}}W\,.
$$
Let $V(W)=\Hom_G(W,V) \otimes_{\mathbb{C}}W$ denote the $W$-isotypical component. Then one clearly has $\det V= \bigotimes_{W\in\hat{G}}\det V(W)$. Let
\begin{equation*}
\det (V,G)=\bigoplus_{W\in\hat{G}}\det V(W)
\end{equation*}
denote the equivariant determinant of $V$. For a $G$-invariant metric on $V$ we define the corresponding equivariant metric on  $\det(V,G)$ as the formal sum
$$
\log|\cdot|^2_{\det(V,G)}:=\sum_{W\in\Irr(G,\mathbb{C})}\log|\cdot|^2_{\det
V(W)}\frac{\chi_W}{\dim_{\mathbb{C}}W}\,,
$$
where $\chi_W$ is the character of $W$.

All this applies as well to the graded representation $V^{\bullet}=H^{\bullet}(M,\mathcal{F})$. The $L^2$-metric on $H^{\bullet}(M,\mathcal{F})$ (viewed as harmonic forms inside $\mathcal{A}^{\bullet}(M,\mathcal{F})$) is a $G$-invariant metric. Let $|\cdot|_{\det(H^{\bullet}(M,\mathcal{F}),G)}$ denote the corresponding equivariant metric on $\det(H^{\bullet}(M,\mathcal{F}),G)$, which we will refer to as the equivariant $L^2$-metric.

Finally, the {\em equivariant Ray-Singer metric} is defined as the formal sum
\begin{equation*}
\|\cdot\|_{\det (H^{\bullet}(M,\mathcal{F}),G)}:=|\cdot|_{\det (H^{\bullet}(M,\mathcal{F}),G)}\cdot\tau\bigl(M,\mathcal{F};\gM,\hF\bigr)\,.
\end{equation*}
We will be interested in the dependence of the equivariant Ray-Singer metric
on the metrics $\gM$ and $\hF$. 

We give ourselves 1-parameter families of $G$-invariant metrics
\begin{enumerate}
\item[(1)] $\varepsilon \mapsto\gM(\varepsilon)$ with $\gM(0)=\gM$
\item[(2)] $\varepsilon \mapsto\hF(\varepsilon)$ with $\hF(0)=\hF$
\end{enumerate}
and we will study the variation of the equivariant Ray-Singer metric
$$
\frac{\partial}{\partial
 \varepsilon}\Bigr\vert_{\varepsilon=0}\log\|\cdot\|_{(\det
H^{\bullet}(M,\mathcal{F}),G)}^2
$$
in each case.

\begin{prop}\label{Vari}
For $\gamma \in G$ one has for the variation of the equivariant Ray-Singer metric:\\
\\
(1) $\varepsilon \mapsto\gM(\varepsilon)$:
\begin{equation*}
 \frac{\partial}{\partial
 \varepsilon}\Bigr\vert_{\varepsilon=0}\log\|\cdot\|_{(\det
H^{\bullet}(M,\mathcal{F}),G)}^2(\gamma)= \LIM_{t\rightarrow 0} \Tr \bigl\{(-1)^{\N}C\gamma\exp(-t\lF))\bigr\}
\end{equation*}
(2) $\varepsilon \mapsto\hF(\varepsilon)$:
\begin{equation*}
 \frac{\partial}{\partial
 \varepsilon}{\Bigr\vert_{\varepsilon=0}}\log\|\cdot\|_{(\det
H^{\bullet}(M,\mathcal{F}),G)}^2(\gamma)= \LIM_{t\rightarrow 0} \Tr \bigl\{(-1)^{N}V\gamma\exp(-t\lF)\bigr\}
\end{equation*}
Here $C=\star^{-1}\dot{\star}$ and $V=(\hF)^{-1}\dot{h}^{\F}$ and
$\LIM_{t\rightarrow 0}f(t)$ denotes the $t^0$-coefficient of the asymptotic
expansion of $f$ as $t\rightarrow 0$ (assuming there exists such).
\end{prop}
\begin{pf}
We use the obvious equivariant generalization of \cite[Prop.~9.38]{BGV},
namely:
\begin{prop}
\label{Var}
Let $\varepsilon\mapsto H(\varepsilon)$ be a 1-parameter family of
$G$-invariant generalized
Laplacians with $H(0)=H$. Let further
\begin{equation*}
\theta(\varepsilon,\gamma,s)=\M\bigl[\;\Tr\bigl\{(-1)^{\N}\N\gamma\exp(-tH(\varepsilon))(1-P_0(\varepsilon))\bigr\}\,\bigr](s)\,.
\end{equation*}
Assume that $\dim \ker H(\varepsilon)$ is constant. Then one has
\begin{equation*}
%\begin{split}
\frac{\partial}{\partial\varepsilon}\Bigr\vert_{\varepsilon=0}\frac{\partial}{\partial
s}\Bigr\vert_{s=0}\theta(\varepsilon,\gamma,s)=-\LIM_{t\rightarrow
0}\Tr\bigl\{(-1)^{\N}\N\gamma\dot{H}H^{-1}\exp(-tH)(1-P_0)\bigr\}\,. 
%\end{split}
\end{equation*}
\end{prop}
The assertion of Proposition \ref{Vari} follows by substituting the formulas
$$
\dot{\Delta}(\F) = -C\dF^*\dF + \dF^* C\dF - \dF C\dF^* + \dF\dF^*C
$$
$$
 \frac{\partial}{\partial
 \varepsilon}{\Bigr\vert_{\varepsilon=0}}\log|\cdot|_{(\det
H^{\bullet}(M,\mathcal{F}),G)}^2(\gamma)=\Tr\bigl\{(-1)^{\N}C\gamma
 P_0\bigr\}   
$$
in case (1) and 
$$
\dot{\Delta}(\F) = -V\dF^*\dF + \dF^* V\dF - \dF V\dF^* + \dF\dF^*V  
$$
$$
 \frac{\partial}{\partial
 \varepsilon}{\Bigr\vert_{\varepsilon=0}}\log|\cdot|_{(\det
H^{\bullet}(M,\mathcal{F}),G)}^2(\gamma) = \Tr\bigl\{(-1)^{\N}V\gamma P_0\bigr\} 
$$
in case (2), whose verification we leave to the reader. 
\end{pf}\\
\\
Therefore Theorem \ref{main} is implied by:
\begin{thm}\label{Asymp}
For $\gamma \in G$ one has:
\begin{xxalignat}{2}
&\text{(1)}\quad\qquad\lim_{t\rightarrow
0}\Tr\bigl\{(-1)^{\N} C\gamma\exp(-t\lF)\bigr\}=-\int\limits_{M^{\gamma}}\theta(\gamma,\mathcal{F},\hF)\wedge\tilde{e}^{\prime}(TM^{\gamma})&&\\
&\text{(2)}\quad\qquad\lim_{t\rightarrow
0}\Tr\bigl\{(-1)^{\N}
V\gamma\exp(-t\lF)\bigr\}=\int\limits_{M^{\gamma}}\Tr[\gamma^{\mathcal{F}} V]e(TM^{\gamma},\nabla^{TM^{\gamma}})&&
\end{xxalignat}
\end{thm}
In the following we will explain the terms, which appear on the right-hand
side of the variation formulas:

Let $or(TM)$ denote the flat line bundle
associated to the orientation cover and $\Pf$ the Pfaffian polynomial. Then
\begin{equation*}
e(TM,\nabla^{TM})=\Pf\Bigl[\,\frac{1}{2\pi}\,R^{TM}\Bigr]
\end{equation*}
is the Euler form of $TM$ associated with $\nabla^{TM}$, where $\nabla^{TM}$ is the Levi-Civita connection on $(M,g^{TM})$ and $R^{TM}$ its curvature. The Euler form is a closed form
and represents via Chern-Weil theory the Euler class $e(TM)\in
H^n(M,or(TM))$.

For a 1-parameter family of Riemannian metrics $\varepsilon\mapsto
g^{TM}(\varepsilon)$ we set
\begin{equation*}
\dot{S}:=\dot{\nabla}^{TM}-\frac{1}{2}\bigl[\nabla^{TM},(g^{TM})^{-1}\dot{g}^{TM}\bigr]\in\mathcal{A}(M,\mathfrak{so}(TM))\,.
\end{equation*}
We define the transgression form
\begin{equation*}
\tilde{e}^{\prime}(TM):=\frac{\partial}{\partial
 b}{\Bigr\vert_{b=0}}\Pf\Bigl[\,\frac{1}{2\pi}\bigl(R^{TM}+b\dot{S}\bigr)\Bigr]\in\mathcal{A}^{n-1}(M,or(TM))
\end{equation*}
and via Chern-Weil theory one obtains the transgression formula
\begin{equation*}
\frac{\partial}{\partial
 \varepsilon}{\Bigr\vert_{\varepsilon=0}}e(TM,\nabla^{TM}(\varepsilon))=d\,\tilde{e}^{\prime}(TM)\in\mathcal{A}^n(M,or(TM))\,.
\end{equation*}
Note that for $\dim M$ odd, $e(TM,\nabla^{TM})$ and $\tilde{e}^{\prime}(TM)$
vanish by the usual properties of $\Pf$. Note further that all this applies as
well to the fixed point set $M^\gamma$, which is a smooth manifold.

Since $\hF$ is not necessarily parallel w.r.t.~$\nF$, we may define a second
flat connection $(\nabla^{\mathcal{F}})^{T}$ on $\F$ by the formula
\begin{equation*}\label{A1}
(\nabla^{\mathcal{F}})^{T}=(\hF)^{-1}\circ\nabla^{\mathcal{F}^*}\circ\hF\,,
\end{equation*}
where $\nabla^{\mathcal{F}^*}$ denotes the connection induced by $\nF$ on $\F^*$ and
$\hF: \F \rightarrow \F^*$ the isomorphism induced by $\hF$. Observe that
$(\nabla^{\mathcal{F}})^{T}=\nF$ if and only if $\nF\hF=0$. As in \cite{BZ1} and \cite{BZ2} we set
\begin{equation*}\label{A2}
\omega(\mathcal{F},\hF):=(\nabla^{\mathcal{F}})^{T}-\nF\in\mathcal{A}^1(M,\End(\mathcal{F}))
\end{equation*}
and for $\gamma \in G$
$$
\theta(\gamma,\mathcal{F},\hF):=\Tr
[\gamma^{\mathcal{F}}\omega(\mathcal{F},\hF)]\in\mathcal{A}^1(M^{\gamma})\,.
$$
Proposition 2.6 in \cite{BZ2} shows that
$\theta(\gamma,\mathcal{F},\hF)\in\mathcal{A}^1(M^{\gamma})$ is closed and that its cohomology class does not depend on $\hF$.

%%%%%%%%%%%%%%%%%%%%%%%%%%%%%%%%%%%
% Proof of the variation formulas %
%%%%%%%%%%%%%%%%%%%%%%%%%%%%%%%%%%%
\section{Proof of the variation formulas}

\subsection{Clifford algebras and exterior algebras}
Let $(V,q)$ be a finite dimensional real vector space equipped with a quadratic
form. Let $C(V,q)$ the associated Clifford algebra, i.e.~the associative
algebra generated by $V$ with the relations $v\cdot w+w\cdot v=-2q(v,w)$ for
$v,w\in V$. The Clifford algebra is a $\mathbb{Z}/2$-graded algebra
(a.k.a.~superalgebra), i.e.~$C(V,q)=C(V,q)^+\oplus C(V,q)^-$ with
$$
C^+(V,q)=\langle v_1\cdot\ldots\cdot v_l : v_i\in V, \, l\text{ even}\rangle
$$
and
$$
C^-(V,q)=\langle v_1\cdot\ldots\cdot v_l : v_i\in V, \, l \text{ odd}\rangle\,.
$$
Recall the filtration $C^{\bullet}(V,q)$ of $C(V,q)$ given by
\begin{equation*}
C^k(V,q)=\langle v_1\cdot\ldots\cdot v_l : v_i\in V,l\leq k\rangle
\end{equation*}
for $k \in \mathbb{Z}$. The associated graded algebra $\Gr^{\bullet}C(V,q)$ is
isomorphic to the exterior algebra $\Lambda^{\bullet}V$ via the symbol map
\begin{equation*}
\sigma : \Gr^{\bullet}C(V,q)\longrightarrow \Lambda^{\bullet}V
\end{equation*}
where the $k$-th symbol is given by
\begin{align*}
\sigma_k: C^k(V,q)/C^{k-1}(V,q) &\longrightarrow \Lambda^kV\\
v_1 \cdot\ldots\cdot v_k + C^{k-1}(V,q) &\longmapsto v_1 \wedge\ldots\wedge
v_k\,. 
\end{align*}
Since as a vector space $C(V,q)$ may be identified with $\Gr^{\bullet}C(V,q)$,
we may also interpret the symbol map as a linear isomorphism
$\sigma:C(V,q)\rightarrow\Lambda V$. In particular $\dim C(V,q)=\dim\Lambda V=2^{\dim V}$.

Let in the following
\begin{align*}
V=&(\mathbb{R}^n,\langle\,\;,\;\rangle)\notag\\
-V=&(\mathbb{R}^n,-\langle\,\;,\;\rangle)\\
(V,-V)=&(\mathbb{R}^n\oplus\mathbb{R}^n,\langle\,\;,\;\rangle\varoplus-\langle\,\;,\;\rangle)\notag
\end{align*}
where $\langle\,\;,\;\rangle$ denotes the standard inner product on $\mathbb{R}^n$,
i.e.~$\langle e_i,e_j\rangle = \delta_{ij}$ for the standard basis $\{e_i\}$.
We will denote $(e_i,0) \in V \oplus V$ again by $e_i$ and $(0,e_i) \in V
\oplus V$ by $\hat{e}_i$.
One has the isomorphism of superalgebras
$$
C^{\pm}(V)\hat{\otimes}C^{\pm}(-V)\longrightarrow C^{\pm}(V,-V)\,.
$$
Furthermore, the tensor product of the symbol maps $\sigma: \Gr^{\bullet} C(V)
\rightarrow \Lambda^{\bullet} V$ and $\hat{\sigma}: \Gr^{\bullet} C(-V)
\rightarrow \Lambda ^{\bullet} V$ yields the symbol map
$$
\sigma\varotimes\hat{\sigma}:\Gr^{\bullet}C(V,-V) \longrightarrow \Lambda^{\bullet}
(V\oplus V)\,,
$$
which we will also denote by $\sigma$.

Using interior multiplication $\iota(e_i):\Lambda^{\bullet}V\rightarrow\Lambda^{\bullet-1}V$
and exterior multiplication
$\varepsilon(e_i):\Lambda^{\bullet}V\rightarrow\Lambda^{\bullet+1}V$ we define
representations of $C(V)$ and $C(-V)$ on the exterior algebra:
\begin{align*}
c\,:\,C^{\pm}(V) &\longrightarrow \End^{\pm}\Lambda V\\
e_i &\longmapsto c(e_i):=\varepsilon(e_i)-\iota(e_i)
\end{align*}
and
\begin{align*}
\hat{c}\,:\,C^{\pm}(-V) &\longrightarrow \End^{\pm}\Lambda V\\
e_i &\longmapsto \hat{c}(e_i):=\varepsilon(e_i)+\iota(e_i)
\end{align*}
The tensor product of these representation yields an isomorphism of
superalgebras
$$
c\varotimes\hat{c}:C^{\pm}(V,-V)   \longrightarrow \End^{\pm}\Lambda V
$$
which we will also denote by $c$.
We obtain a supertrace (i.e.~a linear functional vanishing on
supercommutators) on $C(V,-V)$ by setting
$$
\Str(a)=\Str_{\End\,\Lambda V}[c(a)]
$$
for $a\in C(V,-V)$, where $\Str_{\End\Lambda V}$ is the canonical supertrace
on $\End\Lambda V$.

Let the volume element $\omega \in C(V,-V)$ be defined by
$$
\omega=\pi^{n/2}(-1)^{n(n+1)/2}e_1\cdot\ldots\cdot e_n\cdot \hat{e}_1\cdot\ldots\cdot
\hat{e}_n\,.
$$
The proof of the following lemma is elementary and left to the reader:
\begin{lemma}\label{Superspur}One has:
\begin{enumerate}
\item $\Str$ vanishes on $C^{2n-1}(V,-V)$.
\item $c(\omega)=\pi^{n/2}(-1)^{\N}$, in particular $\Str(\omega)=(4\pi)^{n/2}$.
\end{enumerate}
Here $\N$ denotes the number operator, which multiplies a homogeneous form
with its degree.
\end{lemma}
We will also denote the image of the volume element in $\Lambda (V\oplus V)$ by
$\omega$. For $\alpha \in \Lambda(V\oplus V)$ let $\T\alpha$ be the
coefficient of $\omega$ in $\alpha$. The linear functional $T: \Lambda (V
\oplus V) \rightarrow \mathbb{R}$ is called Berezin trace.
\begin{cor}
For $a \in C(V,V)$ one has $\Str(a)=(4\pi)^{n/2}(\T \circ\,\sigma)(a)$.
\end{cor}
 
\subsection{A Lichnerowicz formula}

In general, neither of the connections $\nF$ and $(\nF)^T$ will preserve the
metric $\hF$. As in \cite{BZ1} we define a third connection
$\nabla^{\mathcal{F},e}=\frac{1}{2}\bigl(\nF+(\nF)^T\bigr)$ on $\F$. This
connection will preserve $\hF$, but it will in general not be flat. 

In the following we will write $\mathcal{E}=\Lambda T^*M\otimes\mathcal{F}$. We
will also denote by $\nabla^{\mathcal{F},e}$ the tensor product connection
$
\nabla^{\Lambda T^*M}\varotimes 1 + 1 \varotimes \nabla^{\mathcal{F},e}
$ 
on $\mathcal{E}$, where $\nabla^{\Lambda T^*M}$ is the connection on $\Lambda T^*M$ induced by $\nabla^{TM}$.
Let $\Delta^{\mathcal{E},e}$ denote the connection Laplacian on $\mathcal{E}$
associated to the connection $\nabla^{\mathcal{E},e}$, i.e.~w.r.t.~a local
ON-frame $\{e_i\}$ one has
\begin{equation*}
\Delta^{\mathcal{E},e}=-\sum_i\Bigl(\bigl(\nabla_{e_i}^{\mathcal{E},e}\bigr)^2-\nabla^{\mathcal{E},e}_{\nabla^{TM}_{e_i}e_i}\Bigr)\,.
\end{equation*}
Since $\nabla^{\mathcal{F},e}$ is a metric connection on $\mathcal{E}$, the
operator $\Delta^{\mathcal{E},e}$ will be formally selfadjoint.

\begin{prop}[\cite{BZ1}: Lichnerowicz formula for $\lF$]\label{Lichnerowicz} One has
$$
\lF = \Delta^{\mathcal{E},e} + E
$$
with $E \in \Gamma(M,\End \mathcal{E})$ which w.r.t.~a local ON-frame $\{e_i\}$ is given by
\begin{align*}
E       =&-\frac{1}{8}\sum_{i,j,k,l}\bigl(R^{TM}(e_i,e_j)e_k,e_l\bigl)c(e_i)c(e_j)\hat{c}(e_k)\hat{c}(e_l)\\
        &-\frac{1}{8}\sum_{i,j}c(e_i)c(e_j)\omega(\mathcal{F},\hF)^2(e_i,e_j)+\frac{1}{8}\sum_{i,j}\hat{c}(e_i)\hat{c}(e_j)\omega(\mathcal{F},\hF)^2(e_i,e_j)\\
        &-\frac{1}{2}\sum_{i,j}c(e_i)\hat{c}(e_j)\Bigl\{\nabla_{e_i}^{T^*M\otimes\End\mathcal{F}}\omega(\mathcal{F},\hF)(e_j)+\frac{1}{2}\omega(\mathcal{F},\hF)^2(e_i,e_j)\Bigr\}\\
        &+\frac{1}{4}\sum_i\bigl(\omega(\mathcal{F},\hF)(e_i)\bigr)^2+\frac{1}{4}r^M\,,
\end{align*}
where $r^M$ denotes the scalar curvature of $(M,g^{TM})$.
\end{prop}
\begin{pf}
We refer the reader to \cite{BZ1}.
\end{pf}

\subsection{Geometry of the frame bundle}

Let us first assume that $M$ is oriented. This assumption will be removed later. Then we can consider the bundle of oriented ON-frames $SO(M)$, i.e.~for $x\in M$ the fiber $SO(M)_x$ consists of all orientation preserving isometries $q:V\rightarrow T_xM$. With the right action $q\mapsto q\cdot h$, $(qh)(v):=q(hv)$, the frame bundle acquires the structure of a $SO(V)$-principal bundle. If $q: M \supset U \rightarrow SO(M)|_U$, $x \mapsto q(x)$ is a local ON-frame and $\{e_i\}$ a (positively oriented) ON-basis of $V$, then $\{qe_i\}$ will be local ON-frame for $TM$, which we will often also denote by $\{e_i\}$ for simplicity. Finally, note that $TM$ may be recovered as associated bundle
$$
TM=SO(M) \times_{SO(V)} V\,.
$$
In the following we will write $Q=SO(M)$ for short. Let $\pi: Q \rightarrow M$ denote the bundle projection.

Let $\omega \in \mathcal{A}^1(Q,\mathfrak{so}(V))$ the 1-form of the connection on $Q$ induced by the Levi-Civita connection on $M$. The 1-form $\omega$ is $SO(V)$-equivariant, i.e.
$$
R_h^*\,\omega = Ad(h)^{-1}\omega
$$
for all $h \in SO(V)$. Further, if the fundamental (vertical) vectorfield $A^Q$ associated with $A\in \mathfrak{so}(V)$ is given by
$$
(A^Q)_q=\left.\frac{d}{dt}\right\vert_{t=0}\,q\cdot\exp(tA)\,,
$$
then one has
$$
\omega(A^Q)=A\,.
$$
A choice of connection on $Q$ yields an $SO(V)$-invariant splitting
$$
TQ=VQ \oplus HQ\,,
$$
where the vertical bundle $VQ$ is given by $VQ=\ker( \pi_*: TQ \rightarrow TM)$ and the horizontal distribution $HQ$ by $HQ=\ker \omega$. For a vectorfield $X$ on $M$ let $X^H$ denote the horiziontal lift of $X$, i.e.~the unique horizontal vectorfield on $Q$ which projects to $X$.

We define the fundamental 1-form $\theta \in \mathcal{A}^1(Q,V)$ by
\begin{equation*} 
\theta_q(X)=q^{-1}(\pi_* X_q) 
\end{equation*} 
for $X \in \Gamma(Q,TQ)$. As the connection 1-form, $\theta$ satisfies an equivariance property, namely $R_h^*\,\theta = h^{-1}\theta$. For $v \in V$, let $v^Q$ denote the fundamental (horizontal) vectorfield associated with $v$, i.e.~the unique horizontal vectorfield on $Q$ which satisfies $\theta(v^Q)=v$.

Let $(\cdot\,,\cdot)_V$ resp.~$(\cdot\,,\cdot)_{\mathfrak{so}(V)}$ denote the inner products on $V$, resp.~on $\mathfrak{so}(V)$. Via the bundle isomorphisms 
$$
Q \times V \rightarrow HQ,\,(q,v) \mapsto v^Q(q)
$$
resp.
$$
Q \times \mathfrak{so}(V) \rightarrow VQ,\,(q,A) \mapsto A^Q(q)
$$
we obtain a Riemannian metric $g^{TQ}$ on $Q$. 

Let $\Omega \in \mathcal{A}^2(Q,\mathfrak{so}(V))$ denote the curvature 2-form 
of $\omega$. Recall that $\Omega$ is $SO(V)$-equivariant, i.e.~$R_h^*\,\Omega = Ad(h)^{-1}\Omega$ for all $h \in SO(V)$, and horizontal, i.e.~$\iota(A^Q)\Omega = 0$ for all $A\in\mathfrak{so}(V)$. For $A \in \mathfrak{so}(V)$ we define $(\Omega,A) \in \mathcal{A}^2(Q)$ by
\begin{equation*}
(\Omega,A)(X,Y)=(\Omega(X,Y),A)_{\mathfrak{so}(V)}
\end{equation*} 
for $X,Y \in \Gamma(Q,TQ)$. Let further denote $\tau:\Lambda^2V\rightarrow \mathfrak{so}(V)$ the unique isomorphism satisfying $(v,\tau(\alpha)w)_V=\alpha(v,w)$ for all $v,w\in V$. Applying this fiberwise, we obtain $\tau\bigl(\Omega,A\bigr)\in\Gamma(Q,\mathfrak{so}(HQ))=C^{\infty}(Q,\mathfrak{so}(V))$.
\begin{lemma}\label{Kommut}One has the following commutator identities:
\begin{enumerate}
\item $[A^Q,v^Q]=\bigl(Av)^Q\,,\;A\in \mathfrak{so}(V)\,,\,v\in V$
\item $[A^Q,B^Q]=[A,B]^Q\,,\;A,B\in\mathfrak{so}(V)$ 
\item   $\omega\bigl([v^Q,w^Q]\bigr)=-\Omega(v^Q,w^Q)\,,\;v,w\in V$\\
        $\theta\bigl([v^Q,w^Q]\bigr)=0$

\item $[A^Q,X^H]=0\,,\;X\in\Gamma(M,TM),A\in\mathfrak{so}(V)$
\item $[X^H,Y^H]=[X,Y]^H-\Omega(X^H,Y^H)^Q\,,\;X,Y\in\Gamma(M,TM)$
\end{enumerate}
\end{lemma}
\begin{pf}
We refer the reader to \cite[Lm.~5.2]{BGV} and \cite[p.~320]{BV}.
\end{pf}\\
\\
This allows us to compute the Levi-Civita connection on $(Q,g^{TQ})$:

\begin{lemma}\label{Levi-Civita}
For the Levi-Civita connection $\nabla^{TQ}$ on $(Q,g^{TQ})$ one has:
\begin{enumerate}
\item $\nabla^{TQ}_{A^Q}B^Q =\frac{1}{2}[A,B]^Q$
\item $\nabla^{TQ}_{v^Q}A^Q =-\frac{1}{2}\tau(\Omega,A)v^Q$\,,\quad $\nabla^{TQ}_{X^H}A^Q =-\frac{1}{2}\tau(\Omega,A)X^H$
\item $\nabla^{TQ}_{A^Q}v^Q = (Av)^Q-\frac{1}{2}\tau(\Omega,A)v^Q$\,,\quad $\nabla^{TQ}_{A^Q}X^H = -\frac{1}{2}\tau(\Omega,A)X^H$
\item $\nabla^{TQ}_{v^Q}w^Q = -\frac{1}{2}\Omega(v^Q,w^Q)^Q$\,,\quad $\nabla^{TQ}_{X^H}Y^H = \bigl(\nabla^{TM}_XY\bigr)^H-\frac{1}{2}\Omega(X^H,Y^H)^Q$
\end{enumerate}
Here $A,B \in \mathfrak{so}(V)$, $v,w \in V$ and $X,Y \in \Gamma(M,TM)$.
\end{lemma}
\begin{pf}
Use the Koszul formula for the Levi-Civita connection together with Lemma \ref{Kommut}.
\end{pf}
\begin{cor}
The trajectories of $v^Q$ and $A^Q$ are geodesics for $v\in V$ and $A \in \mathfrak{so}(V)$. Horizontal lifts of geodesics are geodesics. 
\end{cor}
We record for further reference:
\begin{lemma}\label{Rechts}
For $h \in SO(V)$ one has:
\begin{enumerate}
\item $\exp_q(X)^H\cdot h=\exp_{qh}(X^H)\,,\quad X \in \Gamma(M,TM)$
\item $(\exp_q v^Q)\cdot h=\exp_{qh}\bigl((h^{-1}v)^Q\bigr)\,,\quad v \in V$
\item $\exp_q A^Q\cdot h =q\cdot\exp(A)h\,,\quad A\in\mathfrak{so}(V)$
\end{enumerate}
\end{lemma}
\begin{pf}
Use the fact that $R_h$ is an isometry for $h \in SO(V)$. 
\end{pf}\\
\\
Similarly, in the presence of the group action:
\begin{lemma}\label{Links}
For $\gamma\in G$ one has:
\begin{enumerate}
\item $\gamma \exp_qX^H = \exp_{\gamma q}\bigl((d\gamma)X\bigr)^H\,,\quad X \in \Gamma(M,TM)$
\item $\gamma \exp_qv^Q = \exp_{\gamma q}v^Q\,,\quad v \in V$
\item $\gamma \exp_qA^Q = \exp_{\gamma q}A^Q\,,\quad A\in\mathfrak{so}(V)$
\end{enumerate}
\end{lemma}
\begin{pf}
Use the fact that $L_{\gamma}$ is an isometry for $\gamma \in G$.
\end{pf}\\
\\
Let $q\in SO(M)$. Let $J(q,A): T_qQ\rightarrow T_{q\cdot\exp A}Q$ denote the
differential of $\exp_q$ at $A^Q_q \in T_qQ$. Clearly $J(q,0)=1$. We may view
$J(q,A)$ as endomorphism of $V\oplus \mathfrak{so}(V)$. Moreover, $J(q,A)$
respects this direct sum decomposition, cf.~\cite[Thm.~5.4]{BGV}. As usual,
let $d(q_0,q_1)$ denote the geodesic distance between points $q_0,q_1 \in SO(M)$. 
We define a quadratic form $Q$ on $V\oplus V$ by
\begin{equation*}
Q(A)(v,w)=(v,J(q,A)^{-1}v)-2(v,J(q,A)^{-1}w)+(w,J(q,A)^{-1}w)\,.
\end{equation*}
Clearly one has $Q(0)(v,w)=\|v-w\|^2$.
\begin{lemma}\label{geod_dist}
For $A \in \mathfrak{so}(V)$ small and $v,w \in V$ one has
\begin{align*}
d^2(\exp_qtv^Q,\exp_qtw^Q\cdot\exp A)=\|A\|^2+t^2Q(A)(v,w)+o(t^2)\,.
\end{align*}
\end{lemma}
\begin{pf}
We refer the reader to \cite[Prop.~6.17]{BGV}.
\end{pf}

\subsection{The heat equation}
Let $(M,g^{TM})$ be a closed Riemannian manifold, $(\mathcal{E},h^{\mathcal{E}})$
a hermitian vector bundle over $M$ and $H$ a formally selfadjoint generalized
Laplacian acting on sections of $\mathcal{E}$. For $t>0$, let $(x_0,x_1)
\mapsto k_t(x_0,x_1)$ denote the integral kernel of the heat operator
$\exp(-tH)$. 

One has the well-known asymptotic expansion of the heat kernel, which we will
describe below:
Let us fix points $x_0, x_1 \in M$ and consider the exponential map $\exp_{x_0}:
T_{x_0}M \rightarrow M$. Let $y=\exp_{x_0}^{-1}(x_1)$ and consider the
geodesic $x_s=\exp_{x_0}sy$ connecting $x_0$ to $x_1$.  As usual,
let $d(x_0,x_1)$ denote the geodesic distance between $x_0$ and $x_1$.

\begin{prop}\label{asymp_exp}
There exist $\Phi_i \in \Gamma(M \times M,\mathcal{E}\boxtimes\mathcal{E}^*)$
such that
\begin{align*}
k_t(x_0,x_1) \underset{t\rightarrow0}{\sim} (4\pi t)^{-\dim
M/2}&\exp(-d(x_0,x_1)^2/4t)\Psi_M(d(x_0,x_1)^2)\\
&\cdot j_M(x_0,x_1)^{-1/2}\sum_{i=0}^{\infty}t^i\Phi_i(x_0,x_1)\,,
\end{align*}
where $j_M(x_0, \cdot)$ is the Jacobian of the exponential map at $x_0$ and
$\Psi_M$ a suitable cut-off function.
The coefficients
$x\mapsto\Phi_i(x):=\Phi_i(x_0,x)\in\mathcal{E}_0\otimes\Gamma(X,\mathcal{E}^*)$
are inductively determined by the radial ODE along $x_s$
\begin{equation*}
s\nabla^{\mathcal{E}^*}_{d/ds}\bigl(\Phi_i(x_s) s^i \bigr)=-s^i
j_M^{1/2}H\bigl(j_M^{-1/2}\Phi_{i-1}\bigr)(x_s)
\end{equation*}
with initial condition
\begin{equation*}
\Phi_0(x_0,x_0)=\Id\,.
\end{equation*}
In particular, $\Phi_0(x_0,x_1)=\tau(x_0,x_1)$, the parallel transport along $x_s$. 
\end{prop}
\begin{pf}
We refer the reader to \cite[Thm.~2.30]{BGV} and \cite[p.~329]{BV}.
\end{pf}\\
\\
We may write $\Lambda T^*M$ as an associated bundle
$$
\Lambda T^*M = Q \times_{(SO(V),\Lambda)}\Lambda V
$$
and hence identify sections over $M$ with {\em invariant} sections over $Q$, i.e.
\begin{align*}
\Gamma(M,\Lambda T^*M\otimes\mathcal{F})                &\overset{\cong}{\longrightarrow}
 \bigl(\Lambda V\otimes\Gamma(Q,\pi^*\mathcal{F})\bigr)^{SO(V)}\\
s=\alpha \varotimes \varphi             &\longmapsto \bigl(f_s:q\mapsto(\Lambda q)^{-1}\alpha(x)\varotimes\varphi(x)\bigr)\,.
\end{align*} 
We wish to extend the action of $H:=\lF$ on $\Gamma(M,\Lambda
T^*M\otimes\mathcal{F})$ to the action of a generalized Laplacian $\tilde{H}$ on $\Lambda
V\otimes\Gamma(Q,\pi^*\mathcal{F})$. 
We use the Lichnerowicz formula $\lF = \Delta^{\mathcal{E},e} + E$, cf.~Proposition \ref{Lichnerowicz}, and observe
that the action of $E$ trivially extends. In the following lemma we use  the
{\em Casimir operator} of the representation
$\lambda:\mathfrak{so}(V)\rightarrow \End\Lambda V$, which is given by
$$
\Cas\lambda= \sum_{i<j}(\lambda E_{ij})^2
$$
with $\{E_{ij}\}$ an ON-basis of $\mathfrak{so}(V)$.
\begin{lemma}
One has $\Delta^{\mathcal{E},e}=\bigl(\Delta^{\pi^*\mathcal{F},e}+\Cas\lambda\bigr)\bigr\vert_{\Gamma(M,\mathcal{E})}$.
\end{lemma}
\begin{pf}
For $\{e_i\}$ a local ON-frame one has
\begin{equation*}\begin{split}
\Delta^{\mathcal{E},e}s         &=-\sum_i\Bigl(\bigl(\nabla^{\mathcal{E},e}_{e_i}\bigr)^2-\nabla^{\mathcal{E},e}_{\nabla^{TM}_{e_i}e_i}\Bigr)s\\
                                &=-\sum_i\Bigl(\bigl(\nabla^{\pi^*\mathcal{F},e}_{e_i^H}\bigr)^2-\nabla^{\pi^*\mathcal{F},e}_{(\nabla^{TM}_{e_i}e_i)^H}\Bigr)f_s\,.
\end{split}\end{equation*}
Similarly, since $\{e_i^H,E_{ij}^Q\}$ is a local ON-frame for $TQ$, one has
\begin{align*}
\Delta^{\mathcal{\pi^*\mathcal{F}},e}f_s                &=-\sum_i\Bigl(\bigl(\nabla^{\pi^*\mathcal{F},e}_{e_i^H}\bigr)^2-\nabla^{\pi^*\mathcal{F},e}_{\nabla^{TQ}_{e_i^H}e_i^H}\Bigr)f_s-\sum_{i<j}\Bigl(\bigl(\nabla^{\pi^*\mathcal{F},e}_{E_{ij}^Q}\bigr)^2-\nabla^{\pi^*\mathcal{F},e}_{\nabla^{TQ}_{E_{ij}^Q}E_{ij}^Q}\Bigr)f_s\\
                                                &=-\sum_i\Bigl(\bigl(\nabla^{\pi^*\mathcal{F},e}_{e_i^H}\bigr)^2-\nabla^{\pi^*\mathcal{F},e}_{(\nabla^{TM}_{e_i}e_i)^H}\Bigr)f_s-\sum_{i<j}\bigl(\nabla^{\pi^*\mathcal{F},e}_{E_{ij}^Q}\bigr)^2f_s\\
                                                &=\Delta^{\mathcal{E},e}s-\sum_{i<j}\bigl(\lambda E_{ij}\bigr)^2f_s\,.
\end{align*}
Note that the second line follows from Lemma \ref{Levi-Civita}, i.e. 
$$
\nabla^{TQ}_{E_{ij}^Q}E_{ij}^Q=\frac{1}{2}[E_{ij},E_{ij}]^Q=0 \quad\text{and}\quad \nabla^{TQ}_{e_i^H}e_i^H = (\nabla^{TM}_{e_i}e_i)^H\,.
$$ 
The third line follows from invariance of $f_s$, i.e.
\begin{equation*}
\nabla^{\pi^*\mathcal{F},e}_{E_{ij}^Q}f_s=\left.\frac{d}{dt}\right\vert_{t=0}f_s(q\exp E_{ij}t)=-(\lambda E_{ij})f_s\,.
\end{equation*}
Finally, substitute the definition of $\Cas\lambda$.
\end{pf}\\
\\
We may now set
$$
\tilde{H}=\Delta^{\pi^*\mathcal{F},e}+\Cas\lambda+E\,.
$$
\begin{prop}[Lichnerowicz formula for $\tilde{H}$]\label{Lichnerowicz_lifted} One has
$$
\tilde{H}=\Delta^{\pi^*\mathcal{F},e} + \tilde{E}
$$
with $\tilde{E} \in C(V,-V)\otimes\Gamma(Q,\pi^*\End\mathcal{F})$ which w.r.t.~an ON-basis $\{e_i\}$ of $V$ is given by
\begin{align*}
\tilde{E}_q=    &-\frac{1}{4}\sum_{i,j}c(e_i)c(e_j)\hat{c}(e_i)\hat{c}(e_j)\\
                &-\frac{1}{8}\sum_{i,j,k,l}\bigl(R^{TM}_x(qe_i,qe_j)qe_k,qe_l\bigl)c(e_i)c(e_j)\hat{c}(e_k)\hat{c}(e_l)\\
                &-\frac{1}{8}\sum_{i,j}c(e_i)c(e_j)\omega(\mathcal{F},\hF)^2(qe_i,qe_j)+\frac{1}{8}\sum_{i,j}\hat{c}(e_i)\hat{c}(e_j)\omega(\mathcal{F},\hF)^2(qe_i,qe_j)\\
                &-\frac{1}{2}\sum_{i,j}c(e_i)\hat{c}(e_j)\Bigl\{\nabla_{qe_i}^{T^*M\otimes\End\mathcal{F}}\omega(\mathcal{F},\hF)(qe_i)+\frac{1}{2}\omega(\mathcal{F},\hF)^2(qe_i,qe_j)\Bigr\}\\
                &+\frac{1}{4}\sum_i\bigl(\omega(\mathcal{F},\hF)(qe_i)\bigr)^2+\frac{1}{4}r^M_x-\frac{1}{4}n^2\,.
\end{align*}
\end{prop}
\begin{pf}
An easy calculation yields that
\begin{equation*}
\Cas \lambda = -\frac{1}{4}\sum_{i,j}c(e_i)c(e_j)\hat{c}(e_i)\hat{c}(e_j)-\frac{1}{4}n^2
\end{equation*}
for $\{e_i\}$ an ON-basis of $V$. Now the formula follows with Proposition \ref{Lichnerowicz}.
\end{pf}\\
\\
Let 
$
(q_0,q_1) \mapsto \tilde{k}_t(q_0,q_1) \in \End\Lambda V\otimes\Hom(\mathcal{F}_{x_1},\mathcal{F}_{x_0})
$
be the heat kernel associated with $\tilde{H}$, and 
$
(q_0,q_1) \mapsto k_t(q_0,q_1) \in \End\Lambda V\otimes\Hom(\mathcal{F}_{x_1},\mathcal{F}_{x_0})
$
the lift of the heat kernel associated with $H$ to $Q$. For $t>0$, $q_0,q_1 \in Q$ and $h\in SO(V)$ one has
\begin{equation*}
{k}_t(q_0,q_1)=\int_{SO(V)}\tilde{k}_t(q_0,q_1h)\Lambda h^{-1} dh\,,
\end{equation*}
hence for $t>0$ and $q\in\pi^{-1}(x)$
\begin{equation*}
\Str[k_t(x,x)]  = \int_{SO(V)}\Str[\tilde{k}_t(q,qh)\Lambda h^{-1}]dh\,.   
\end{equation*} 
We fix $q=q_0 \in Q$ and consider for $A \in
\mathfrak{so}(V)$ the asymptotic expansion of
$\tilde{k}_t(q_0,q_0\exp(A))$. We write $\tilde{\Phi}_i(A)$ instead of
$\tilde{\Phi}_i(q_0,q_0\exp(A))$ for short.
\begin{prop}\label{Transport}
For the coefficients $\tilde{\Phi}_i(A)$, $A \in \mathfrak{so}(V)$, one has:
\begin{enumerate}
\item $\tilde{\Phi}_i(A) \in C^{4i}(V,-V)\otimes \End \mathcal{F}_x$.
\item The sum of the highest symbols for $A=0$ is given by 
\begin{align*}
\sum_{i=0}^{n/2}\bigl(\sigma_{4i}\tilde{\Phi}_i\bigr)(0)=\exp\Bigl(\frac{1}{8}\sum_{i,j,k,l}&\bigl(R_x^{TM}(qe_i,qe_j)qe_k,qe_l\bigl)e_i\wedge
e_j\wedge \hat{e}_k\wedge \hat{e}_l\\
        &+ \frac{1}{4}\sum_{i,j}e_i\wedge
e_j\wedge\hat{e}_i\wedge \hat{e}_j\Bigr)\,.
\end{align*}
\end{enumerate}
\end{prop}
\begin{pf}
1.~We look at the radial ODE determining the coefficients $\tilde{\Phi}_i$
    along the geodesic $q_s:=\exp_{q_0}(sY)$, $Y\in T_{q_0}Q$, cf.~Proposition
\ref{asymp_exp}:
\begin{equation*}
s\nabla^{\pi^*\mathcal{F},e}_{d/ds}\bigl(s^i\tilde{\Phi}_i(q_s)\bigr)=-s^ij_Q^{1/2}\tilde{H}\bigl(j_Q^{-1/2}\tilde{\Phi}_{i-1}\bigr)(q_s)\,.
\end{equation*}
We may divide by $s$ to obtain
\begin{equation}\label{grosernull}
\nabla^{\pi^*\mathcal{F},e}_{d/ds}\bigl(s^i\tilde{\Phi}_i(q_s)\bigr)=-s^{i-1}j_Q^{1/2}\tilde{H}\bigl(j_Q^{-1/2}\tilde{\Phi}_{i-1}\bigr)(q_s)\,.
\end{equation}
Clearly $\tilde{\Phi}_0(q_1)=\tau(q_0,q_1)\in
C^0(V,-V)\otimes\Hom(\mathcal{F}_{x_1},\mathcal{F}_{x_0})$. Assuming by
induction that $\tilde{\Phi}_{i-1}\in
C^{4i-4}(V,-V)\otimes\mathcal{F}_{x_0}\otimes\Gamma(Q,\pi^*\mathcal{F}^*)$ we
obtain using the Lichnerowicz formula
$\tilde{H}=\Delta^{\pi^*\mathcal{F},e}+\tilde{E}$ with $\tilde{E} \in
C^4(V,-V)\otimes\Gamma(Q,\pi^*\End\mathcal{F})$ that the right-hand side of (\ref{grosernull}) lies in
$C^{4i}(V,-V)\otimes\Hom(\mathcal{F}_{x_s},\mathcal{F}_{x_0})$, hence so by (\ref{grosernull}) the coefficient $\tilde{\Phi}_i(q_1)$.

2.~To calculate $\sum_{i=0}^{n/2}\bigl(\sigma_{4i}\tilde{\Phi}_i\bigr)(0)$, we
specialize our considerations to vertical geodesics $q_s=q_0\exp(sA)$. Since
the connection $\nabla^{\pi^*\mathcal{F},e}$ is trivial in fiber direction, (\ref{grosernull}) becomes very simple:
\begin{equation}\label{grosernullII}
d/ds\bigl(s^i\tilde{\Phi}_i(sA)\bigr)=-s^{i-1}j_Q^{1/2}\tilde{H}\bigl(j_Q^{-1/2}\tilde{\Phi}_{i-1}\bigr)(sA)\,.
\end{equation}
For the highest symbols we obtain:
\begin{equation}\label{hSymb}
d/ds\Bigl(s^i\bigl(\sigma_{4i}\tilde{\Phi}_i\bigr)(sA)\Bigr)				=-s^{i-1}\bigl(\sigma_{4}\tilde{E}\bigr)\wedge\bigl(\sigma_{4i-4}\tilde{\Phi}_{i-1}\bigr)(sA)\,.
\end{equation}
We set
$$
f_i(s)	=\bigl(\sigma_{4i}\tilde{\Phi}_i\bigr)(sA)
$$
and 
$$
F(s)	=\sum_{i=0}^{n/2} s^i\,f_i(s)
$$
such that using (\ref{hSymb}) we get:
\begin{equation}\label{Transp}
d/ds\,F(s)=-(\sigma_{4}\tilde{E})(sA)\wedge F(s)\,.
\end{equation}
Observe that $F(0)=\bigl(\sigma_{0}\tilde{\Phi}_i\bigr)(0)=1\varotimes\Id$ and
$F(1)=\sum_{i=0}^{n/2}\bigl(\sigma_{4i}\tilde{\Phi}_i\bigr)(A)$. 
Recall that $\bigl(\sigma_{4}\tilde{E}\bigr)(q)$ is given by
\begin{align*}
\bigl(\sigma_{4}\tilde{E}\bigr)(q_0) =
	&- \frac{1}{8}\sum_{i,j,k,l}\bigl(R_{x_0}^{TM}(q_0e_i,q_0e_j)q_0e_k,q_0e_l\bigl)e_i\wedge e_j
\wedge \hat{e}_k \wedge\hat{e}_l\\
	&- \frac{1}{4}\sum_{i,j}e_i\wedge e_j \wedge
\hat{e}_i \wedge \hat{e}_j\,. 
\end{align*}
The curvature term is equivariant w.r.t.~$SO(V)$, hence
\begin{align*}
\bigl(\sigma_{4}\tilde{E}\bigr)(sA) =
	&- \frac{1}{8}\exp(-s\lambda A)\sum_{i,j,k,l}\bigl(R_{x_0}^{TM}(q_0e_i,q_0e_j)q_0e_k,q_0e_l\bigl)e_i\wedge e_j\wedge \hat{e}_k \wedge\hat{e}_l\\
	&- \frac{1}{4}\sum_{i,j}e_i\wedge e_j \wedge
\hat{e}_i \wedge \hat{e}_j\,, 
\end{align*}
and we get in (\ref{Transp}):
\begin{align*}
d/ds\,F(s)=	&\frac{1}{8}\Bigl(\exp(-s\lambda
A)\sum_{i,j,k,l}\bigl(R_{x_0}^{TM}(q_0e_i,q_0e_j)q_0e_k,q_0e_l\bigl)e_i\wedge
e_j\wedge \hat{e}_k \wedge\hat{e}_l\Bigr)\\
& \wedge F(s)+\frac{1}{4}\Bigl(\sum_{i,j}e_i\wedge e_j \wedge
\hat{e}_i \wedge \hat{e}_j\Bigr) \wedge F(s)
\end{align*}
We deform this system of ODE into a system with constant coefficients.
For $t \in [0,1]$ we consider:
\begin{align*}
d/ds\,F_t(s)=	&\frac{1}{8}\Bigl(\exp(-st\lambda
A)\sum_{i,j,k,l}\bigl(R_{x_0}^{TM}(q_0e_i,q_0e_j)q_0e_k,q_0e_l\bigl)e_i\wedge
e_j\wedge \hat{e}_k \wedge\hat{e}_l\Bigr)\\
& \wedge F_t(s)+\frac{1}{4}\Bigl(\sum_{i,j}e_i\wedge e_j \wedge
\hat{e}_i \wedge \hat{e}_j\Bigr) \wedge F_t(s)\,.
\end{align*}
For $t=0$ we can explicitly solve this:
\begin{align*}
F_0(1)=		\exp\Bigl\{&\frac{1}{8}\Bigl(\sum_{i,j,k,l}\bigl(R_{x_0}^{TM}(q_0e_i,q_0e_j)q_0e_k,q_0e_l\bigl)e_i\wedge e_j\wedge \hat{e}_k \wedge\hat{e}_l\Bigl)\\
		&+\frac{1}{4}\Bigl(\sum_{i,j}e_i\wedge e_j \wedge
\hat{e}_i \wedge \hat{e}_j\Bigr)\Bigl\}(1\varotimes\Id)\,.
\end{align*}
Using continuous dependence of the solution on the coefficients of the ODE and
the continuity of the $\tilde{\Phi}_i$ we get
\begin{equation*}
\sum_{i=0}^{n/2}\bigl(\sigma_{4i}\tilde{\Phi}_i\bigr)(0)
=\lim_{t\rightarrow0}\sum_{i=0}^{n/2}\bigl(\sigma_{4i}\tilde{\Phi}_i\bigr)(At)
=\lim_{t\rightarrow0}F_t(1)
=F_0(1)\,.
\end{equation*} 
This finishes the proof.
\end{pf}\\
\\
In contrast to $V=(\hF)^{-1}\dot{\hF}$, the endomorphism $C=\star^{-1}\dot{\star}$ contains Clifford variables, more precisely one has:
\begin{lemma}
\label{Hodge}
The endomorphism $C=\star^{-1}\dot{\star}$ is given in terms of Clifford variables by
\begin{equation*}
C_q=-\frac{1}{2}\sum_{i,j}\bigl((g^{TM})^{-1}\dot{g}^{TM}e_i,e_j\bigr)_xc(e_i)\hat{c}(e_j) \in C^{2}(V,-V)\,.
\end{equation*}
In particular one has
\begin{equation*}
\bigl(\sigma_{2}C\bigr)(q)=-\frac{1}{2}\sum_{i,j}\bigl((g^{TM})^{-1}\dot{g}^{TM}e_i,e_j\bigr)_x
e_i\wedge \hat{e}_j \in \Lambda^{2}(V\oplus V)\,.
\end{equation*}
\end{lemma}
\begin{pf}
We refer the reader to \cite[Thm.~4.15]{BZ1}.
\end{pf}\\
\\
This fact has to be taken care of in the evaluation of $\lim_{t\rightarrow0}\Str\{C\gamma\exp(-tH)\}$. To facilitate the computations, J.M.~Bismut and W.~Zhang introduce an extra even Clifford variable $\sigma$, cf.~\cite{BZ1}, which will also turn out to be useful in our approach, cf.~Proposition \ref{orient}:

Let $\mathbb{R}\langle\sigma\rangle:=\mathbb{R}[\sigma]\bigl/(\sigma^2-1)$. We consider the trivial $\mathbb{Z}/2$-grading on $\mathbb{R}\langle\sigma\rangle$, i.e.~$\sigma$ is even.
If $W^{\pm}$ is a real  $\mathbb{Z}/2$-graded vector space (a.k.a.~superspace), then $W^{\pm}\otimes\mathbb{R}\langle\sigma\rangle$ is a $\mathbb{R}\langle\sigma\rangle$-module. One has
\begin{align*}
\End_{\mathbb{R}\langle\sigma\rangle}(W^{\pm}\otimes\mathbb{R}\langle\sigma\rangle)&\overset{\cong}{\longrightarrow}\End_{\mathbb{R}}W^{\pm}\otimes\mathbb{R}\langle\sigma\rangle\\
A+\sigma B &\longmapsto A\varotimes 1 + B \varotimes \sigma\notag
\end{align*}
and we can extend the trace by $\mathbb{R}\langle\sigma\rangle$-linearity:
\begin{align*}
\Str_{\mathbb{R}\langle\sigma\rangle}:\End_{\mathbb{R}\langle\sigma\rangle}(W^{\pm}\otimes\mathbb{R}\langle\sigma\rangle)&\rightarrow\mathbb{R}\langle\sigma\rangle\\
A+\sigma B &\mapsto \Str[A] + \sigma\Str[B]\,.\notag
\end{align*}
In this situation we denote:
\begin{align*}
\Str^1[A+\sigma B]		&:=\Str[A]\,,\\
\Str^{\sigma}[A+\sigma B]	&:=\Str[B]\,.
\end{align*}
We define
\begin{equation*}
H_{odd}:=-\frac{1}{2}\sum_{i,j}c(e_i)\hat{c}(e_j)\Bigl\{\bigl(\nabla_{e_i}\omega\bigr)(e_j)+\frac{1}{2}\omega^2(e_i,e_j)\Bigr\}
\end{equation*}
and $H_{ev}:=H-H_{odd}$. Then $H_{ev}+\sigma H_{odd}$ is a generalized Laplacian on the vector bundle $\mathcal{E}\otimes\mathbb{R}\langle\sigma\rangle$ with heat kernel $(x_0,x_1) \mapsto l_t(x_0,x_1)$.

We fix $q=q_0 \in Q$ and consider for $A \in
\mathfrak{so}(V)$ the asymptotic expansion of
$\tilde{l}_t(q_0,q_0\exp(A))$. We write $\tilde{\Phi}_i(A)$ instead of
$\tilde{\Phi}_i(q_0,q_0\exp(A))$ for short.
\begin{prop}\label{sigmaTransport}
For the coefficients $\tilde{\Phi}_i(A)$, $A \in \mathfrak{so}(V)$, one has:
\begin{enumerate}
\item $\tilde{\Phi}_i(A)\in
C(V,-V)\otimes\End\mathcal{F}_x\otimes\mathbb{R}\langle\sigma\rangle$,
i.e.~$\tilde{\Phi}_i(A)=\tilde{\Phi}^1_i(A)+\sigma\tilde{\Phi}_i^{\sigma}(A)$
with $\tilde{\Phi}_i^{1}, \tilde{\Phi}_i^{\sigma} \in C(V,-V)\otimes\End\mathcal{F}_x$.

\item $\tilde{\Phi}^{\sigma}_i(A) \in C^{4i-2}(V,-V)\otimes \End \mathcal{F}_x$.

\item The sum of the highest symbols for $A=0$ is given by 
\begin{align*}
\sum_{i=1}^{n/2}\bigl(\sigma_{4i-2}\tilde{\Phi}^{\sigma}_i\bigr)(0)=&\exp\Bigl(\frac{1}{8}\sum_{i,j,k,l}\bigl(R_x^{TM}(qe_i,qe_j)qe_k,qe_l\bigl)e_i\wedge
e_j\wedge \hat{e}_k\wedge \hat{e}_l\\
        &\qquad\qquad\quad+\frac{1}{4}\sum_{i,j}e_i\wedge
e_j\wedge\hat{e}_i\wedge \hat{e}_j\Bigr)\\
        &\wedge\Bigl(\frac{1}{2}\sum_{i,j}e_i\wedge \hat{e}_j\Bigl\{\bigl(\nabla_{qe_i}\omega\bigr)(qe_j)+\frac{1}{2}\omega^2(qe_i,qe_j)\Bigr\}\Bigr)\,.
\end{align*}
\end{enumerate}
\end{prop}
\begin{pf}
We proceed as in the proof of Proposition \ref{Transport}:

1.~The initial condition of the radial ODE yields 
$$
\tilde{\Phi}_0(q_1)=\tau(q_0,q_1)\in C(V,-V)\otimes\Hom(\mathcal{F}_{x_1},\mathcal{F}_{x_0})\otimes\mathbb{R}\langle\sigma\rangle\,.
$$
It is $\tilde{E}_{ev}+\sigma\tilde{E}_{odd} \in
C(V,-V)\otimes\Gamma(Q,\pi^*\End\mathcal{F})\otimes\mathbb{R}\langle\sigma\rangle$.
The assertion follows by induction.

2.~The initial condition of the radial ODE yields $\tilde{\Phi}_0^1(q_1)=\tau(q_0,q_1)$ and $\tilde{\Phi}_0^{\sigma}(q_1)=0$. Observe that
\begin{align*}
\bigl(\tilde{E}_{ev}+\sigma\tilde{E}_{odd}\bigl)\bigl(\tilde{\Phi}_{i-1}^1+\sigma\tilde{\Phi}_{i-1}^{\sigma}\bigr)
=&\tilde{E}_{ev}\tilde{\Phi}_{i-1}^1+\tilde{E}_{odd}\tilde{\Phi}_{i-1}^{\sigma}\\
&+\sigma\bigl(\tilde{E}_{odd}\tilde{\Phi}_{i-1}^1+\tilde{E}_{ev}\tilde{\Phi}_{i-1}^{\sigma}\bigr)\,.
\end{align*}
Furthermore, recall that $\tilde{E}_{ev}\in C^4(V,-V)\otimes\Gamma(Q,\pi^*\End\mathcal{F})$ and that $\tilde{E}_{odd}\in C^2(V,-V)$ $\otimes\Gamma(Q,\pi^*\End\mathcal{F})$. Then the assertion follows by induction.

3.~We obtain for the highest symbols:
\begin{align*}
d/ds\Bigl(s^i\bigl(\sigma_{4i}\tilde{\Phi}_i^1\bigr)(sA)\Bigr)				&=-s^{i-1}\bigl(\sigma_{4}\tilde{E}\bigr)\wedge\bigl(\sigma_{4i-4}\tilde{\Phi}^{1}_{i-1}\bigr)(sA)\,,\\
d/ds\Bigl(s^i\bigl(\sigma_{4i-2}\tilde{\Phi}_i^{\sigma}\bigr)(sA)\Bigr)	
	&=-s^{i-1}\Bigl(\bigl(\sigma_{4}\tilde{E}\bigr)(sA)\wedge\bigl(\sigma_{4i-6}\tilde{\Phi}_{i-1}^{\sigma}\bigl)(sA)\\
	&\qquad\qquad+\bigl(\sigma_{2}\tilde{E}\bigr)(sA)\wedge\bigl(\sigma_{4i-4}\tilde{\Phi}_{i-1}^1\bigl)(sA)\Bigl)\,.
\end{align*}
We set
\begin{align*}
f_i^{1}(s) &:=\bigl(\sigma_{4i}\tilde{\Phi}_i^{1}\bigr)(sA)\,,\,i\geq0\,,\\
f_i^{\sigma}(s) &:=\bigl(\sigma_{4i-2}\tilde{\Phi}_i^{\sigma}\bigr)(sA)\,,\,i\geq1
\end{align*}
and
\begin{align*}
F^1(s)		&:=\sum_{i=0}^{n/2}s^if_i^1(s)\,,\\
F^{\sigma}(s)	&:=\sum_{i=1}^{n/2}s^if_i^{\sigma}(s)\,.	
\end{align*}
Then we have
\begin{alignat}{2}
F^1(0)&=1\varotimes\Id&\quad\text{and}\quad F^1(1)&=\sum_{i=0}^{n/2}\bigl(\sigma_{4i}\tilde{\Phi}_i^1\bigr)(A)\,,\\
F^{\sigma}(0)&=0&\quad\text{and}\quad F^{\sigma}(1)&=\sum_{i=1}^{n/2}\bigl(\sigma_{4i-2}\tilde{\Phi}_i^{\sigma}\bigr)(A)
\end{alignat}
and we obtain
\begin{align}
d/ds F^1(s)		&=-\bigl(\sigma_{4}\tilde{E}\bigr)(sA)\wedge F^1(s)\label{1}\,,\\
d/ds F^{\sigma}(s)	&=-\bigl(\sigma_{4}\tilde{E}\bigr)(sA)\wedge
F^{\sigma}(s)-\bigl(\sigma_{2}\tilde{E}\bigr)(sA)\wedge F^1(s)\,.\label{sigma} 
\end{align}
It is
\begin{align*}
\bigl(\sigma_{4}\tilde{E}\bigr)(sA) =
	&- \frac{1}{8}\exp(-s\lambda A)\sum_{i,j,k,l}\bigl(R_x^{TM}(qe_i,qe_j)qe_k,qe_l\bigl)e_i\wedge e_j\wedge \hat{e}_k \wedge\hat{e}_l\\
	&- \frac{1}{4}\sum_{i,j}e_i\wedge e_j \wedge
\hat{e}_i \wedge \hat{e}_j 
\end{align*}
and
\begin{align*}
\bigl(\sigma_{2}\tilde{E}\bigr)(sA)=-\frac{1}{2}\exp(-s\lambda A)\sum_{i,j}e_i\wedge\hat{e}_j\Bigl\{\bigl(\nabla_{qe_i}\omega\bigr)(qe_j)+\frac{1}{2}\omega^2(qe_i,qe_j)\Bigl\}\,.
\end{align*}
As in the proof of Proposition \ref{Transport} we set $A=0$ and get in (\ref{1}), (\ref{sigma}):
\begin{align}
d/ds F^1(s)	=&\frac{1}{8}\Bigl(\sum_{i,j,k,l}\bigl(R_x^{TM}(qe_i,qe_j)qe_k,qe_l\bigl)e_i\wedge e_j\wedge \hat{e}_k \wedge\hat{e}_l\Bigr)\wedge F^1(s)\label{1.2}\\
		&+\frac{1}{4}\Bigl(\sum_{i,j}e_i\wedge e_j \wedge
\hat{e}_i \wedge \hat{e}_j\Bigr) \wedge F^1(s)\notag\,;
\end{align}
\begin{align}
d/ds F^{\sigma}(s)
=&\frac{1}{8}\Bigl(\sum_{i,j,k,l}\bigl(R_x^{TM}(qe_i,qe_j)qe_k,qe_l\bigl)e_i\wedge e_j\wedge \hat{e}_k \wedge\hat{e}_l\Bigr)\wedge F^{\sigma}(s)\label{sigma.2}\\
		&+\frac{1}{4}\Bigl(\sum_{i,j}e_i\wedge e_j \wedge
\hat{e}_i \wedge \hat{e}_j\Bigr) \wedge F^{\sigma}(s)\notag\\
&+\frac{1}{2}\Bigl(\sum_{i,j}e_i\wedge\hat{e}_j\Bigl\{\bigl(\nabla_{qe_i}\omega\bigr)(qe_j)+\frac{1}{2}\omega^2(qe_i,qe_j)\Bigr\}\Bigr)\wedge
F^1(s)\notag\,.
\end{align}
From Proposition \ref{Transport} we know the solution for (\ref{1.2}), which yields the inhomogeneity in (\ref{sigma.2}):
\begin{align*}
F^1(s)=		\exp\Bigl\{&\frac{1}{8}\Bigl(\sum_{i,j,k,l}\bigl(R_x^{TM}(qe_i,qe_j)qe_k,qe_l\bigl)e_i\wedge e_j\wedge \hat{e}_k \wedge\hat{e}_l\Bigl)s\\
		&+\frac{1}{4}\Bigl(\sum_{i,j}e_i\wedge e_j \wedge
\hat{e}_i \wedge \hat{e}_j\Bigr)s\Bigl\}(1\varotimes\Id)\,.
\end{align*}
For (\ref{sigma.2}) we obtain:
\begin{align*}
F^{\sigma}(1)=\exp\Bigl(&\frac{1}{8}\sum_{i,j,k,l}\bigl(R_x^{TM}(qe_i,qe_j)qe_k,qe_l\bigl)e_i\wedge e_j\wedge \hat{e}_k \wedge\hat{e}_l\\
		&+\frac{1}{4}\sum_{i,j}e_i\wedge e_j \wedge
\hat{e}_i \wedge \hat{e}_j\Bigr)\\
&\cdot\Bigl(\frac{1}{2}\sum_{i,j}e_i\wedge\hat{e}_j\Bigl\{\bigl(\nabla_{qe_i}\omega\bigr)(qe_j)+\frac{1}{2}\omega^2(qe_i,qe_j)\Bigr\}\Bigr)(1\varotimes\Id)\,.
\end{align*}
This finishes the proof.
\end{pf}\\
\\
In the presence of the group action, let $(x_0,x_1) \mapsto
k_t(\gamma,x_0,x_1)$ denote the integral kernel of the operator
$\gamma\exp(-tH)$. For $t>0$ and $x_0,x_1 \in M$ one has
$$
k_t(\gamma,x_0,x_1) = \gamma^{\mathcal{E}} k_t(\gamma^{-1}x_0,x_1)\,,
$$ 
hence for $t>0$ 
$$
\Str\{\gamma\exp(-tH)\}=\int_M\Str[k_t(\gamma,x)]dx\,,
$$
where we write $k_t(\gamma,x)$ for $k_t(\gamma,x,x)$.

Since $M$ is a closed manifold, the fixed point set $M^\gamma$ of the isometry
$\gamma$ is a disjoint union of finitely many components $M^\gamma_i$, which
are compact submanifolds without boundary, possibly of varying dimension, cf.~\cite[Thm.~5.1]{Kob}. Since our
calculations are local, we may assume that the fixed point set consists of a
single component.

The tangent bundle of $M$ decomposes orthogonally over $M^\gamma$ as
$$
TM\vert_{M^{\gamma}}=TM^{\gamma}\oplus\mathcal{N}\,,
$$
where $\mathcal{N}$ denotes the normal bundle of $M^\gamma$ in $M$. Clearly
$TM^{\gamma}$ is precisely the eigenbundle of $d\gamma$ corresponding to the
eigenvalue $1$. Let $n_0=\dim M^\gamma$ and $n_1=n-n_0$. We write $V=V_0+V_1$
with $V_0\cong\mathbb{R}^{n_0}$ and $V_1\cong\mathbb{R}^{n_1}$.

Let $\phi \in C_c^{\infty}(M)$ be a function which is equal to $1$ on
$M^{\gamma}$ and vanishes outside a tubular neighbourhood $U$ of $M^{\gamma}$. For $V\in\Gamma(M,\End\mathcal{F})$ we set
\begin{equation*}
I(t,\gamma,\phi,x)=\int_{V_1}\Str\bigl[Vk_t(\gamma,\exp_xqv)\bigr]\phi(\exp_xqv)dv\,.
\end{equation*}
\begin{prop}\label{Vasymp}
For $x\in M^{\gamma}$ there are $\Phi_l \in
C_c^{\infty}(\mathfrak{so}(V),C(V,-V)\otimes\End\mathcal{F}_x)$ such that
\begin{align*}
I(t,\gamma,\phi,x)\underset{t\rightarrow0}{\sim}(4\pi t)^{(n_1-\dim
Q)/2}\sum_{l=0}^{\infty}
t^l\int_{\mathfrak{so}(V)}&\exp(-\|A\|^2/4t)\\
&\cdot\Str\bigl[(\Lambda\tilde{\gamma})\Phi_l(A)\exp(-\lambda
A)\bigr]dA
\end{align*}
with $\Phi_l(A)\in C^{4l}(V,-V)\otimes\End\mathcal{F}_x$ and highest symbol
\begin{align*}
\sigma_{4l}\Phi_l(A)=V\gamma^{\pi^*\mathcal{F}}\bigl(\sigma_{4l}\tilde{\Phi}_l\bigr)(A)\Psi_{\mathfrak{so}(V)}(A)j_{\mathfrak{so}(V)}(A)\det
Q_1(A,\gamma)^{-1/2}\,.
\end{align*}
\end{prop}
\begin{pf}
Let $q\in SO(M)$ such that $q(V_0)\subset T_xM^{\gamma}$ and
$q(V_1)\subset\mathcal{N}_x$. Let further in the following denote
$x_v=\exp_xqv$ and $q_v=\exp_qv^Q$. We then have $\pi(q_v)=x_v$ and
\begin{align*}
I(t,\gamma,\phi,x) &= \int_{V_1}\Str\bigl[Vk_t(\gamma,x_v)\bigr]\phi(x_v)dv\\
&=\int_{V_1}\Str\bigl[V\gamma^{\pi^*\mathcal{F}}k_t(\gamma^{-1}q_v,q_v)\bigr]\phi(q_v)dv\\
&=\int_{V_1}\int_{SO(V)}\Str\bigl[V\gamma^{\pi^*\mathcal{F}}\tilde{k}_t(q_v,\gamma
q_vh)\Lambda
h^{-1}\bigr]\phi(q_v)dhdv\,.
\end{align*}
Let $\tilde{\gamma}=\tilde{\gamma}(q)\in SO(V)$ be determined by the requirement that $\gamma q=q\tilde{\gamma}$.
With Lemma \ref{Links} and Lemma \ref{Rechts} we get
$\gamma\exp_qv^Q=\exp_{\gamma
q}v^Q=\exp_q(\tilde{\gamma}v)^Q\tilde{\gamma}$, hence $\gamma
q_v=q_{\tilde{\gamma}v}\tilde{\gamma}$. Upon substituting
$\tilde{\gamma}h\mapsto h$ we obtain
\begin{align*}
I(t,\gamma,\phi,x)
&=\int_{V_1}\int_{SO(V)}\Str\bigl[V\gamma^{\pi^*\mathcal{F}}\tilde{k}_t(q_v,q_{\tilde{\gamma}v}h)\Lambda(h^{-1}\tilde{\gamma})\bigr]\phi(q_v)dhdv\,.
\end{align*}
Asymptotically, as $t\rightarrow0$, we may replace the integration over $SO(V)$
by an integration over the Lie algebra $\mathfrak{so}(V)$ and substitute the
asymptotic expansion for $\tilde{k}_t$, cf.~Proposition \ref{asymp_exp}. With $\Lambda\bigl(\exp(-A)\tilde{\gamma}\bigr)=\exp(-\lambda
A)(\Lambda\tilde{\gamma})$ we get
\begin{align*}
\Str\bigl[V\gamma^{\pi^*\mathcal{F}}\tilde{\Phi}_k(A)\Lambda\bigl(\exp(-A)\tilde{\gamma}\bigr)\bigr]
=&\Str\bigl[(\Lambda\tilde{\gamma})V\gamma^{\pi^*\mathcal{F}}\tilde{\Phi}_k(A)\exp(-\lambda
A)\bigr]
\end{align*}
and further
\begin{align*}
I(t,\gamma,\phi,x)
&\underset{t\rightarrow0}{\sim}
\int_{V_1}\int_{\mathfrak{so}(V)}\Str\bigl[(\Lambda\tilde{\gamma})V\gamma^{\pi^*\mathcal{F}}\tilde{k}_t(q_v,q_{\tilde{\gamma}v}\exp(-\lambda
A)\bigr]\phi(q_v)\\
&\qquad\qquad\qquad\qquad\qquad\qquad\qquad\qquad\quad\cdot\Psi_{\mathfrak{so}(V)}(A)j_{\mathfrak{so}(V)}(A)dAdv\\
&\underset{t\rightarrow0}{\sim}(4\pi t)^{-\dim
Q/2}\sum_{k=0}^{\infty}t^k\int_{V_1}\int_{\mathfrak{so}(V)}\exp\bigl(-d(q_v,q_{\tilde{\gamma}v}\exp
A)^2/4t\bigr)\\
&\qquad\cdot\Str\bigl[(\Lambda
\tilde{\gamma})V\gamma^{\pi^*\mathcal{F}}\tilde{\Phi}_k(A,v)\exp(-\lambda A)\bigr]\phi(v)\Psi(A,v)j(A,v)dAdv\,.
\end{align*}
Let $h(A,v)=d^2\bigl(q_v,q_{\tilde{\gamma}v}\exp A\bigr)$. From Lemma \ref{geod_dist} we obtain that $v=0$ is a critical point of $v\mapsto
h_A(v):=h(A,v)$, and that $Q_1(A):=\frac{1}{2}\operatorname{Hess}\bigl\vert_{v=0}h_A(v)$ is
positive definite. The Morse lemma ensures the existence of local coordinates
$w(v)=F_A(v)$ about $0\in V$ with $w(0)=0$ such that
$h_A(w)=\|A\|^2+\|w\|^2$. We further have $|\det dF_A(0)|=\det Q_1(A)^{1/2}$
and $Q_1(0)(v)=\|(1-\tilde{\gamma})v\|^2$.

Let $m(a,w)$ be the Jacobian of the coordinate change,
i.e.~$m(A,w)dw=dv$. Then it follows that $m(A,0)=\det Q_1(A)^{-1/2}$. We may
interchange the order of integration and get
\begin{align*}
I(t,\gamma,\phi,x)      &\underset{t\rightarrow0}{\sim}(4\pi t)^{-\dim
Q/2}\sum_{k=0}^{\infty}t^k\int_{\mathfrak{so}(V)}\int_{V_1}\exp\bigl(-(\|A\|^2+\|w\|^2)/4t\bigr)\\
&\qquad\qquad\qquad\qquad\qquad\qquad\quad\cdot\Str\bigl[(\Lambda\tilde{\gamma})V\gamma^{\pi^*\mathcal{F}}\tilde{\Phi}_k(A,w)\exp(-\lambda
A)\bigr]\\
&\qquad\qquad\qquad\qquad\qquad\qquad\quad\cdot\phi(w)\Psi(A,w)j(A,w)m(A,w)\,dwdA\\
&\underset{t\rightarrow0}{\sim}(4\pi t)^{-\dim
Q/2}\sum_{k=0}^{\infty}t^k\int_{\mathfrak{so}(V)}\int_{V_1}\exp\bigl(-(\|A\|^2+\|w\|^2)/4t\bigr)\\
&\qquad\qquad\qquad\qquad\qquad\qquad\quad\cdot\Str\bigl[(\Lambda\tilde{\gamma})f_k(A,w)\exp(-\lambda
A)\bigr]\,dwdA\,,
\end{align*}
where $f_k(A,w)$ is given by
\begin{equation*}
f_k(A,w):=V\gamma^{\pi^*\mathcal{F}}\tilde{\Phi}_k(A,w)\phi(w)\Psi(A,w)j(A,w)m(A,w)\,.
\end{equation*}
We use the asymptotic expansion
\begin{align*}
(4\pi t)^{-n_1/2}\int_{V_1}\exp(-\|w\|^2/4t)f_k(A,w)dw\underset{t\rightarrow0}{\sim}\,\sum_{i=0}^{\infty}f_{k,i}(A)t^{i} 
\end{align*}
with
\begin{align*}
f_{k,0}(A)=f_k(A,0)=V\gamma^{\pi^*\mathcal{F}}\tilde{\Phi}_k(A)\Psi_{\mathfrak{so}(V)}(A)j_{\mathfrak{so}(V)}(A)\det
Q_1(A,\gamma)^{-1/2}
\end{align*}
and we set
\begin{equation*}
\Phi_l(A)=\sum_{i=0}^lf_{l-i,i}(A)\,.
\end{equation*}
From Proposition \ref{Transport} we deduce that $\Phi_l(A) \in
C^{4l}(V,-V)\otimes\End{F}_x$ and that the highest symbol is given by
\begin{align*}
\sigma_{4l}\Phi_l(A) &= \sigma_{4l}f_l(A,0)\\
&= V\gamma^{\pi^*\mathcal{F}}\bigl(\sigma_{4l}\tilde{\Phi}_l\bigr)(A)\Psi_{\mathfrak{so}(V)}(A)j_{\mathfrak{so}(V)}(A)\det
Q_1(A,\gamma)^{-1/2}\,.
\end{align*}
Altogether we obtain
\begin{align*}
I(t,\gamma,\phi,x)      &\underset{t\rightarrow0}{\sim}(4\pi t)^{-\dim
Q/2}\sum_{k=0}^{\infty}\sum_{i=0}^{\infty}t^{k+i}\int_{\mathfrak{so}(V)}\exp\bigl(-(\|A\|^2)/4t\bigr)\\
&\qquad\qquad\qquad\qquad\qquad\qquad\quad\cdot\Str\bigl[(\Lambda\tilde{\gamma})f_{k,i}\exp(-\lambda
A)\bigr]dA\\
&\underset{t\rightarrow0}{\sim}(4\pi t)^{-\dim
Q/2}\sum_{l=0}^{\infty}t^l\int_{\mathfrak{so}(V)}\exp\bigl(-(\|A\|^2)/4t\bigr)\\
&\qquad\qquad\qquad\qquad\qquad\qquad\cdot\Str\bigl[(\Lambda\tilde{\gamma})\Phi_l(A,w)\exp(-\lambda
A)\bigr]dA\,,
\end{align*}
which finishes the proof.
\end{pf}\\
\\
Similarly, for $C\in\Gamma(M,C^2(M,-M))$ we set
\begin{equation*}
I^{\sigma}(t,\gamma,\phi,x)=\int_{V_1}\Str^{\sigma}\bigl[Cl_t(\gamma,\exp_xqv)\bigr]\phi(\exp_xqv)dv
\end{equation*}
and we obtain in the same way as in the proof of Proposition \ref{Vasymp}
(using Proposition \ref{sigmaTransport} instead of Proposition \ref{Transport}):
\begin{prop}\label{Casymp}
For $x\in M^{\gamma}$ there are $\Phi_l^{\sigma} \in
C_c^{\infty}(\mathfrak{so}(V),C(V,-V)\otimes\End\mathcal{F}_x)$ such that
\begin{align*}
I^{\sigma}(t,\gamma,\phi,x)\underset{t\rightarrow0}{\sim}(4\pi t)^{(n_1-\dim
Q)/2}\sum_{l=0}^{\infty}
t^l\int_{\mathfrak{so}(V)}&\exp(-\|A\|^2/4t)\\
&\cdot\Str\bigl[(\Lambda\tilde{\gamma})\Phi_l^{\sigma}(A)\exp(-\lambda
A)\bigr]dA
\end{align*}
with $\Phi_l^{\sigma}(A)\in C^{4l}(V,-V)\otimes\End\mathcal{F}_x$ and highest symbol
\begin{align*}
\sigma_{4l}\Phi_l^{\sigma}(A)=\gamma^{\pi^*\mathcal{F}}\bigl(\sigma_2
C\bigr)(A)\wedge&\bigl(\sigma_{4l-2}\tilde{\Phi}_l^{\sigma}\bigr)(A)\\
&\cdot\Psi_{\mathfrak{so}(V)}(A)j_{\mathfrak{so}(V)}(A)\det
Q_1(A,\gamma)^{-1/2}\,.
\end{align*}
\end{prop}

\subsection{Asymptotics of Gaussian integrals}
Let us call a multiindex $\alpha=(\alpha_1,\ldots,\alpha_N) \in \mathbb{N}_0^N$
even, if all $\alpha_i$ are even numbers. We call a multiindex $\alpha$
odd, if at least one $\alpha_i$ is an odd number.
\begin{lemma}\label{gaussInt}
Let $\alpha=(\alpha_1,\ldots,\alpha_N) \in \mathbb{N}_0^N$. Then one has:
\begin{equation*}
(4\pi t)^{-N/2}\int\limits_{\mathbb{R}^N} \exp(-\|x\|^2/4t)x^{\alpha}dx=\left\{ 
\begin{array}{cl}
        t^{|\alpha|/2}\frac{\alpha!}{(\alpha/2)!} & \text{if $\alpha$ is even}\\
        0 & \text{if $\alpha$ is odd}
\end{array}\right.
\end{equation*}
\end{lemma}
\begin{pf}
One has
$$
(4\pi t)^{-N/2}\int\limits_{\mathbb{R}^N}\exp(-\|x\|^2/4t)x^{\alpha}dx=\frac{\partial^{\alpha}}{\partial b^{\alpha}}\Bigr\vert_{b=0}\exp(t\|b\|^2)\,.
$$
To evaluate this expression, we look at the power series expansion
\begin{equation*}
\exp(t\|b\|^2)=\sum_{k=0}^{\infty}\frac{t^k(\sum_{i=1}^Nb_i^2)^k}{k!}\,.
\end{equation*}
The coefficient of $b^{\alpha}=b_1^{\alpha_1}\cdot\ldots\cdot b_N^{\alpha_N}$ in
this expansion is $0$, if $\alpha$ is odd, and it is $t^{|\alpha|/2}/(\alpha/2)!$, if $\alpha$ is even.
\end{pf}
\begin{cor}\label{koeff}
For $i \in \mathbb{N}_0$ one has:
\begin{xalignat*}{2}
\lim_{t\rightarrow0}(4\pi t)^{-N/2}t^{-i}\int\limits_{\mathbb{R}^N}\exp(-\|x\|^2/4t)x^{\alpha}dx=
\left\{\begin{array}{cl}
        0&\text{$\alpha$ odd}\\
        0&\text{$\alpha$ even, $i<|\alpha|/2$}\\
        \frac{\alpha!}{(\alpha/2)!}&\text{$\alpha$ even, $i=|\alpha|/2$}\\
        \infty&\text{$\alpha$ even, $i>|\alpha|/2$}
\end{array}\right.
\end{xalignat*}
\end{cor}
Let $\varphi \in C_c^{\infty}(\mathbb{R}^N)$. The Taylor expansion about $x=0$ is an asymptotic expansion as $x\rightarrow0$:
\begin{equation*}
\varphi(x) \underset{x\rightarrow0}{\sim} \sum_{k=0}^{\infty}\sum_{|\alpha|=k}\varphi_{\alpha}x^{\alpha}\,,
\end{equation*}
where the coefficient $\varphi_{\alpha}$ is given by
\begin{equation*}
\varphi_{\alpha}=\frac{1}{\alpha!}\Bigl(\frac{\partial^{\alpha}}{\partial x^{\alpha}}\varphi\Bigr)(0)\,.
\end{equation*}
From this we obtain, using Lemma \ref{gaussInt} term by term, an asymptotic expansion as $t\rightarrow 0$, cf.~\cite[p.~73]{BGV}:
\begin{equation*}\begin{split}\label{asymp}
        &(4\pi t)^{-N/2}t^{-i}\int\limits_{\mathbb{R}^N}\exp(-\|x\|^2/4t)\varphi(x)dx\\
         \underset{t\rightarrow0}{\sim}&\,(4\pi
         t)^{-N/2}\sum_{k=0}^{\infty}\sum_{|\alpha|=k}\varphi_{\alpha}t^{-i}\int\limits_{\mathbb{R}^N}\exp(-\|x\|^2/4t)x^{\alpha}dx\\
\underset{t\rightarrow0}{\sim}&\,\sum_{k=0}^{\infty}\sum_{|\alpha|=2k \atop \alpha\,\text{even}}\varphi_{\alpha}\frac{\alpha!}{(\alpha/2)!}t^{k-i}
\end{split}\end{equation*}

\begin{lemma}\label{asymp2}
If the coefficients $\varphi_{\alpha}$ vanish for all even $\alpha$ with $|\alpha|<2i$, then one has
\begin{equation*}
\lim_{t\rightarrow0}\,(4\pi
t)^{-N/2}t^{-i}\int\limits_{\mathbb{R}^N}\exp(-\|x\|^2/4t)\varphi(x)dx =
\sum_{|\alpha|=2i \atop \alpha\,\text{even}}\varphi_{\alpha}\frac{\alpha!}{(\alpha/2)!}\,.
\end{equation*}
\end{lemma}
\begin{pf}
We apply Corollary \ref{koeff} to the above expansion term by term. The conditions on the $\varphi_{\alpha}$ ensure that the singular terms vanish. The formula for the limit is also obtained from Corollary \ref{koeff} as the sum of the nonvanishing terms.
\end{pf}\\
\\
Let us now set $N=n(n-1)/2$ and identify $\mathbb{R}^N$ with $\mathfrak{so}(V)$. Let further $\phi \in C_c^{\infty}(\mathfrak{so}(V),C(V,-V))$ be given. We will further assume that $\phi(A)$ is an element in $C^{2i,2i}(V,-V)$ for all $A \in \mathfrak{so}(V)$. 
We claim that under these conditions the limit
\begin{equation*}
\lim_{t\rightarrow0}\,(4\pi
t)^{-\dim
Q/2}t^{i}\int_{\mathfrak{so}(V)}\exp(-\|A\|^2/4t)\Str\bigl[\phi(A)\exp(-\lambda A)\bigr]dA
\end{equation*}
exists.

We introduce some further notation. The space of polynomial functions on $\mathfrak{so}(V)$ may by fixing the basis $\{E_{ij}\}$ be identified with the ring of polynomials in $n(n-1)/2$ variables, which we will denote by $\mathbb{R}[A_{ij}]$. Let $\E=\sum_{i<j}E_{ij}\in\mathfrak{so}(V)$. Evaluation of a polynomial $p \in \mathbb{R}[A_{ij}]$ at $\E$ yields the sum of the coefficients:
\begin{equation*}
p(\E)=\sum_{\alpha}p_{\alpha}A^{\alpha}(A=\E)=\sum_{\alpha}p_{\alpha}\,.
\end{equation*}

We define a formal power series $Q(A)$ with coefficients in the Clifford algebra by
\begin{equation*}
Q(A)=\exp(-\lambda A) \in
\mathbb{R}\llbracket A_{ij} \rrbracket\otimes C(V,-V)\,.
\end{equation*}
Recall that
$
-\lambda A =
\frac{1}{2}\sum_{i<j}A_{ij}(c(e_i)c(e_j)-\hat{c}(e_i)\hat{c}(e_j))$
and hence by taking the highest symbol
$
\sigma_2(-\lambda A)=\frac{1}{2}\sum_{i<j}A_{ij}(e_i\wedge e_j-\hat{e}_i\wedge\hat{e}_j)
$.

Let the formal power series $P(A)$ with coefficients in the even part of the exterior algebra be defined by
\begin{equation*}
P(A)=\exp\bigl(\sigma_2(-\lambda A)\bigr)\in
\mathbb{R}\llbracket A_{ij} \rrbracket\otimes\Lambda^{2\bullet}(V\oplus V)\,. 
\end{equation*}
One clearly has the identity
\begin{equation*}
P(A)=\sum_{\alpha}(\sigma_{2|\alpha|}Q_{\alpha})A^{\alpha}=\exp\Bigr(\frac{1}{2}\sum_{i<j}A_{ij}\bigl(e_i\wedge
e_j-\hat{e}_i\wedge\hat{e}_j\bigr)\Bigr)\,.
\end{equation*}
\begin{lemma}
One has the identity
\begin{equation*}
\sum_{\alpha\,\text{even}}P_{\alpha}\frac{\alpha!}{(\alpha/2)!}A^{\alpha}=\exp\Bigl(-\frac{1}{2}\sum_{i<j}A_{ij}^2e_i\wedge
e_j\wedge \hat{e}_i\wedge \hat{e}_j\Bigr)\in\mathbb{R}\llbracket A_{ij}
\rrbracket\otimes\Lambda^{2\bullet}(V\oplus V)\,.
\end{equation*}
\end{lemma}
\begin{pf}
The assertion follows by comparing coefficients. Details are left to the reader.
\end{pf}\\
\\
We have to look at a slightly more general situation. Let $F$ be a finite dimensional, trivially graded vector space. We extend the trace, resp.~the Berezin trace, in the obvious way
\begin{align*}
\Str\varotimes\Tr:C(V,-V)\otimes\End F\longrightarrow\mathbb{R}\\
\T  \varotimes\Tr:\Lambda(V\oplus V)\otimes\End F\longrightarrow\mathbb{R}
\end{align*}
and we will in the following also denote these maps by $\Str$, resp.~by $\T$. Similarly, we will also denote the map
\begin{align*}
\sigma\varotimes\Id:C(V,-V)\otimes\End F\longrightarrow \Lambda(V\oplus
V)\otimes\End F\,.
\end{align*}
by $\sigma$ in the following.
\begin{prop}\label{asymp4}
For $\phi \in C_c^{\infty}(\mathfrak{so}(V),C^{4i}(V,-V)\otimes \End F)$ one has
\begin{equation*}\begin{split}
&\lim_{t\rightarrow0}\,(4\pi t)^{-\dim
Q/2}t^{i}\int_{\mathfrak{so}(V)}\exp(-\|A\|^2/4t)\Str\bigl[\phi(A)\exp(-\lambda
A)\bigr]dA\\
=&\T\Bigl(\exp\bigl(-\frac{1}{2}\sum_{i<j}e_i\wedge
e_j\wedge \hat{e}_i\wedge \hat{e}_j\bigr)\bigl(\sigma_{4i}\phi\bigr)(0)\Bigr)\,.
\end{split}\end{equation*}
\end{prop}
\begin{pf}
We look at the Taylor expansion of
$\varphi(A)=\Str[\phi(A)\exp(\lambda A)]$ about $A=0$
\begin{equation*}
\varphi(A) \underset{A\rightarrow0}{\sim}\,
\sum_{k=0}^{\infty}\sum_{|\alpha|=k}\varphi_{\alpha}A^{\alpha}\in\mathbb{R}\llbracket
A_{ij}\rrbracket
\end{equation*}
and we wish to apply Lemma \ref{asymp2}. Therefore we have to show that the coefficient $\varphi_{\alpha}$ vanishes as long as $\alpha$ is even with $|\alpha|<n-2i$. But this is clear since $\phi(A)\frac{1}{k!}(-\lambda A)^k\in C^{2(2i+k)}(V,-V)\otimes \End F$,
such that
$\Str\bigl[\phi(A)\frac{1}{k!}(-\lambda A)^k\bigr]$
vanishes for $k<n-2i$, cf.~Lemma \ref{Superspur}.
Now Lemma \ref{asymp2} gives us the following expression for the limit:
\begin{equation*}\begin{split}
&\lim_{t\rightarrow0}\,(4\pi t)^{-\dim
Q/2}t^{i}\int_{\mathfrak{so}(V)}\exp(-\|A\|^2/4t)\Str\bigl[\phi(A)\exp(-\lambda
A)\bigr]dA\\
=&(4\pi)^{-n/2}\sum_{|\alpha|=n-2i \atop \alpha\,\text{even}}\varphi_{\alpha}\frac{\alpha!}{(\alpha/2)!}\,.
\end{split}\end{equation*}
Observe that
\begin{align*}
&(4\pi)^{-n/2}\sum_{|\alpha|=n-2i \atop
\alpha\,\text{even}}\varphi_{\alpha}\frac{\alpha!}{(\alpha/2)!}
=(4\pi)^{-n/2}\Str\Bigl[\phi(0)\sum_{|\alpha|=n-2i \atop
\alpha\,\text{even}}Q_{\alpha}\frac{\alpha!}{(\alpha/2)!}\,\Bigr]
\end{align*}
since all other monomials in the $A_{ij}$ of total degree $n-2i$ do not have trace. To see this use that $\phi(A)\in C^{4i}(V,-V)\otimes \End F$ for $A \in \mathfrak{so}(V)$ and that $Q_{\alpha}\in C^{2|\alpha|}(V,-V)\otimes \End F$.
Continuing with the calculation, we get:
\begin{align*}
&(4\pi)^{-n/2}\Str\Bigl[\phi(0)\sum_{|\alpha|=n-2i \atop
\alpha\,\text{even}}Q_{\alpha}\frac{\alpha!}{(\alpha/2)!}\,\Bigr]\\
=&\bigl(\T\circ\,\sigma\bigr)\Bigl(\phi(0)\sum_{|\alpha|=n-2i \atop
\alpha\,\text{even}}Q_{\alpha}\frac{\alpha!}{(\alpha/2)!}\,\Bigr)\\
=&\T\Bigl(\bigl(\sigma_{4i}\phi\bigr)(0)\sum_{|\alpha|=n-2i \atop
\alpha\,\text{even}}\sigma_{2|\alpha|}Q_{\alpha}\frac{\alpha!}{(\alpha/2)!}\,\Bigr)\\
=&\T\Bigl(\bigl(\sigma_{4i}\phi\bigr)(0)\sum_{|\alpha|=n-2i \atop
\alpha\,\text{even}}P_{\alpha}\frac{\alpha!}{(\alpha/2)!}A^{\alpha}\,\Bigr)(A=\E)\\
=&\T\Bigl(\bigl(\sigma_{4i}\phi\bigr)(0)\exp\bigl(-\frac{1}{2}\sum_{i<j}A_{ij}^2e_i\wedge
e_j\wedge \hat{e}_i\wedge \hat{e}_j\bigr)\Bigr)(A=\E)\\
=&\T\Bigl(\exp\bigl(-\frac{1}{2}\sum_{i<j}e_i\wedge
e_j\wedge \hat{e}_i\wedge \hat{e}_j\bigr)\bigl(\sigma_{4i}\phi\bigr)(0)\Bigr)
\end{align*}
In the second last line, we may substitute the full power series, since the
monomials of degree $\neq n-2i$ do not contribute to the Berezin trace. In the
last line, we use that $\Lambda^{2\bullet}(V\oplus V)$ is commutative.
\end{pf}

\begin{cor}\label{PotenzreiheII}
For a formal power series $\Phi(t,A)=\sum_{i=0}^{\infty}t^i\Phi_i(A)$ with coefficients $\Phi_i \in C_c^{\infty}(\mathfrak{so}(V),C^{4i}(V,-V)\otimes\End F)$ one has\begin{equation*}\begin{split}
&\lim_{t\rightarrow0}\,(4\pi t)^{-\dim
Q/2}\sum_{i=0}^{\infty}t^{i}\int_{\mathfrak{so}(V)}\exp(-\|A\|^2/4t)\Str\bigl[\Phi_i(A)\exp(-\lambda
A)\bigr]dA\\
=&\T\Bigl(\exp\bigl(-\frac{1}{2}\sum_{i<j}e_i\wedge
e_j\wedge \hat{e}_i\wedge \hat{e}_j\bigr)\sum_{i=0}^{n/2}\bigl(\sigma_{4i}\Phi_i\bigr(0)\Bigr)\,.
\end{split}\end{equation*}
\end{cor}
\begin{pf}
Use Proposition \ref{asymp4} term by term.
\end{pf}

\subsection{Evaluation of the asymptotic terms}
With $V=V_0 \oplus V_1$ as above we have
\begin{align*}
C(V,-V)&=C(V_0,-V_0)\hat{\otimes}C(V_1,-V_1)\\
\Lambda (V\oplus V)&=\Lambda(V_0\oplus V_0)\hat{\otimes}\Lambda(V_1\oplus V_1)
\end{align*}
and for the volume elements $\omega=\omega_0\cdot\omega_1$, resp.~$\omega=\omega_0\wedge\omega_1$. We write $\Str_{0}$ and $\Str_{1}$, resp.~$\T_{0}$ and $\T_{1}$, for the trace, resp.~the Berezin trace, on the corresponding spaces. Let $p_0$ denote the orthogonal projection $V\rightarrow V_0$. 

We define:
\begin{align*}
\sigma_{k}^0:C^{k,l}(V,-V)      &\longrightarrow \Lambda^{k}(V_0\oplus V_0)\\
                        a       &\longmapsto \bigl((\Lambda p_0)\circ\sigma_{k}\bigr)a\notag 
\end{align*}
Note that monomials $e_{i_1\ldots i_k}\cdot\hat{e}_{j_1\ldots j_l}$ are killed by $\sigma_{k}^0$ if at least one index is larger than $n_0$.
\begin{lemma}
\label{Hallo}
For a formal power series $\Phi(t,A)=\sum_{i=0}^{\infty}t^i\Phi_i(A)$ with coefficients $\Phi_i \in C^{\infty}(\mathfrak{so}(V),C^{4i}(V,-V)\otimes\End F)$ and $a_1\in C(V_1,-V_1)$ one has 
\begin{equation*}\begin{split}
&\lim_{t\rightarrow0}\,(4\pi t)^{(n_1-\dim
Q)/2}\sum_{i=0}^{\infty}t^{i}\int_{\mathfrak{so}(V)}\exp(-\|A\|^2/4t)\Str\bigl[a_1\Phi_i(A)\exp(-\lambda
A)\bigr]dA\\
=&\Str_{1}(a_1)\T_{0}\Bigl(\exp\Bigl\{-\frac{1}{4}\sum_{i,j=1}^{n_0}e_i\wedge
e_j\wedge \hat{e}_i\wedge \hat{e}_j\Bigr\}\sum_{i=0}^{n_0/2}(\sigma^0_{4i}\Phi_i)(0)\Bigr)\,.
\end{split}\end{equation*}
\end{lemma}
\begin{pf}
Since $C(V_1,-V_1)$ is contained in $C^{2n_1}(V,-V)$, it follows that $a_1\Phi_i(A)$ is an element in $C^{2n_1+4i}(V,-V)$. In particular we have
\begin{equation*}
\sigma_{2n_1+4i}\bigl(a_1\Phi_i(A)\bigr)=\sigma_{2n_1}a_1\wedge\sigma_{4i}\Phi(A)
\end{equation*}
for $A\in\mathfrak{so}(V)$. Hence we may apply Corollary \ref{PotenzreiheII} to the formal power series $\Phi^{\prime}(t,A)=\sum_{k=n_1/2}^{\infty}t^ka_1\Phi_{k-n_1/2}(A)$ and we obtain
\begin{align*}
&\lim_{t\rightarrow0}(4\pi)^{n_1/2}(4\pi t)^{-\dim Q/2}\sum_{k=n_1/2}^{\infty}t^k
\int_{\mathfrak{so}(V)}\exp(-\|A\|^2/4t)\\
&\qquad\qquad\qquad\qquad\qquad\qquad\qquad\qquad\qquad\cdot\Str\bigl[a_1\Phi_{k-n_1/2}(A)\exp(-\lambda A)\bigr]dA\\
=&(4\pi)^{n_1/2}\T\Bigl(\exp\Bigl\{-\frac{1}{4}\sum_{i,j=1}^{n/2}e_i\wedge
e_j\wedge\hat{e}_i\wedge\hat{e}_j\Bigr\}\wedge\bigl(\sigma_{2n_1}a_1\bigr)\\
&\qquad\qquad\qquad\qquad\qquad\qquad\qquad\qquad\qquad\;\;\wedge\sum_{k=n_1/2}^{n}\bigl(\sigma_{4k-2n_1}\Phi_{k-n_1/2}\bigr)(0)\Bigr)\\
=&\Str_1(a_1)\T\Bigl(\exp\Bigl\{-\frac{1}{4}\sum_{i,j=1}^{n/2}e_i\wedge
e_j\wedge\hat{e}_i\wedge\hat{e}_j\Bigr\}\wedge\omega_1\wedge\sum_{i=0}^{n_0/2}\bigl(\sigma_{4i}\Phi_{i}\bigr)(0)\Bigr)\,.
\end{align*} 
In the last line we have used that $\sigma_{2n_1}a_1=(4\pi)^{-n_1/2}\Str_1(a_1)\omega_1$. Since $\omega_1$ annulates monomials, which contain $e_i$ or $\hat{e}_i$ with $i>n_0$, the assertion follows from the definition of $\sigma^0_{4i}$.
\end{pf}
\begin{prop}\label{I}
Let $x \in M^{\gamma}$. Then the limit $I(\gamma,x)=\lim_{t\rightarrow0}I(t,\gamma,\phi,x)$ exists and one has
\begin{align*}
I(\gamma,x)=&\Tr[\gamma^{\mathcal{F}}_{x}V_x]\\
&\cdot\T_0\Bigl(\exp\Bigl\{\frac{1}{8}\sum_{i,j,k,l=1}^{n_0}\bigl(R_{x}^{TM^{\gamma}}(qe_i,qe_j)qe_k,qe_l\bigl)e_i\wedge
e_j\wedge\hat{e}_k\wedge\hat{e}_l\Bigr\}\Bigr)\,.
\end{align*}
\end{prop}
\begin{pf}
We wish to show that the asymptotic expansion in Proposition \ref{Vasymp} converges as $t\rightarrow 0$. To that end, we first observe that $\tilde{\gamma}$ acts as the identity on $V_0$, which in turn yields $\Lambda\tilde{\gamma} \in C(V_1,-V_1)$.
Now we apply Lemma \ref{Hallo} with $a_1=\Lambda \tilde{\gamma}$ and
$\Phi_i(A)$ as given by Proposition \ref{Vasymp}. Besides convergence, we get for the limit 
\begin{align*}
I(\gamma,x)=\Str_1(\Lambda\tilde{\gamma})\T_0\Bigl(\exp\Bigl\{-\frac{1}{4}\sum_{i,j=1}^{n_0}e_i\wedge
e_j\wedge\hat{e}_i\wedge\hat{e}_j\Bigr\}\sum_{i=0}^{n_0/2}\bigl(\sigma_{4i}^{0}\Phi_i\bigr)(0)\Bigr)\,.
\end{align*}
Proposition \ref{Vasymp} also yields
\begin{align*}
\bigl(\sigma_{4i}^0\Phi_i\bigr)(0)=V_x\gamma_x^{\pi^*\mathcal{F}}\bigl(\sigma_{4i}^0\tilde{\Phi}_i\bigr)(0)\det
Q_1(0,\gamma)^{-1/2}\,.
\end{align*}
Next, Proposition \ref{Transport} gives us for the sum of the highest symbols
\begin{align*}
\sum_{i=0}^{n_0/2}\bigl(\sigma_{4i}^0\tilde{\Phi}_i\bigr)(0)=&\exp\Bigl\{\frac{1}{4}\sum_{i,j=1}^{n_0}e_i\wedge
e_j\wedge\hat{e}_i\wedge{e}_j\Bigr\}\\
&\wedge\exp\Bigl\{\frac{1}{8}\sum_{i,j,k,l=1}^{n_0}\bigl(R_x^{TM^{\gamma}}(qe_i,qe_j)qe_k,qe_l\bigr)e_i\wedge
e_j\wedge\hat{e}_k\wedge{e}_l\Bigr\}\,.
\end{align*}
Altogether we obtain
\begin{align*}
I(\gamma,x)=&\det Q_1(0,\gamma)^{-1/2}\Str\bigl[\Lambda
\tilde{\gamma}\bigr]\Tr\bigl[V_x\gamma^{\mathcal{F}}_x\bigr]\\
&\cdot\T_0\Bigl(\exp\Bigl\{\frac{1}{8}\sum_{i,j,k,l=1}^{n_0}\bigl(R_{x}^{TM^{\gamma}}(qe_i,qe_j)qe_k,qe_l\bigl)e_i\wedge
e_j\wedge\hat{e}_k\wedge\hat{e}_l\Bigr\}\Bigr)\,.
\end{align*}
Observe now that on the one hand $\Str_{\Lambda V_1}[\Lambda\tilde{\gamma}]=\det_{V_1}(1-\tilde{\gamma})>0$, since $1$ is not an eigenvalue of $\tilde{\gamma}\vert_{V_1}$. On the other hand we have already seen that $\det_{V_1}
Q_1(0,\gamma)=\det_{V_1}\bigl((1-\tilde{\gamma})^T(1-\tilde{\gamma})\bigr)=\det_{V_1}(1-\tilde{\gamma})^2$. Hence theses terms cancel and we get the result.
\end{pf}\\
\\
We will use Lemma \ref{Hodge} in the proof of the following proposition, which corresponds to Proposition \ref{I}: 
\begin{prop}
Let $x \in M^{\gamma}$. Then the limit $I^{\sigma}(\gamma,x)=\lim_{t\rightarrow0}I^{\sigma}(t,\gamma,\phi,x)$ exists and one has 
\begin{align*}
I^{\sigma}(\gamma,x)=-\T_0\Bigl(\,&\frac{1}{2}\Bigl\{\sum_{i,j=1}^{n_0}\bigl((g^{TM^{\gamma}})^{-1}\dot{g}^{TM^{\gamma}}qe_i,qe_j\bigr)_x
e_i\wedge \hat{e}_j\Bigr\}\\
        &\wedge\exp\Bigl\{\frac{1}{8}\sum_{i,j,k,l=1}^{n_0}\bigl(R_x^{TM^{\gamma}}(qe_i,qe_j)qe_k,qe_l\bigl)e_i\wedge
e_j\wedge \hat{e}_k \wedge\hat{e}_l\bigr)\Bigr\}\\
        &\wedge\frac{1}{2}\Bigl\{\sum_{i,j=1}^{n_0}e_i\wedge\hat{e}_j\Tr\bigl[\gamma_x^{\mathcal{F}}\bigl(\nabla_{qe_i}\omega\bigr)(qe_j)\bigr]\Bigr\}\,\Bigr)\,. 
\end{align*}
\end{prop}
\begin{pf} 
We proceed as in the proof of Proposition \ref{I}. Using Proposition \ref{Casymp} (instead of Proposition \ref{Vasymp}) we get
 \begin{align*}
I^{\sigma}(\gamma,x)=\Str_1(\Lambda\tilde{\gamma})\T_0\Bigl(\exp\Bigl\{-\frac{1}{4}\sum_{i,j=1}^{n_0}e_i\wedge
e_j\wedge\hat{e}_i\wedge\hat{e}_j\Bigr\}\sum_{i=0}^{n_0/2}\bigl(\sigma_{4i}^{0}\Phi_i^{\sigma}\bigr)(0)\Bigr)
\end{align*}
and
\begin{align*}
\bigl(\sigma_{4i}^0\Phi_i^{\sigma}\bigr)(0)=\gamma_x^{\pi^*\mathcal{F}}\bigl(\sigma_2^0
C\bigr)(0)\wedge\bigl(\sigma_{4i-2}^0\tilde{\Phi}_i\bigr)(0)\det
Q_1(0,\gamma)^{-1/2}\,.
\end{align*}
Next, Proposition \ref{sigmaTransport} (instead of Proposition \ref{Transport}) gives us
\begin{align*}
\sum_{i=0}^{n_0/2}\bigl(\sigma_{4i-2}^0\tilde{\Phi}_i^{\sigma}\bigr)(0)=&\exp\Bigl\{\frac{1}{4}\sum_{i,j=1}^{n_0}e_i\wedge
e_j\wedge\hat{e}_i\wedge{e}_j\Bigr\}\\
&\wedge\exp\Bigl\{\frac{1}{8}\sum_{i,j,k,l=1}^{n_0}\bigl(R_x^{TM^{\gamma}}(qe_i,qe_j)qe_k,qe_l\bigr)e_i\wedge
e_j\wedge\hat{e}_k\wedge{e}_l\Bigr\}\\
&\wedge\Bigl\{\frac{1}{2}\sum_{i,j=1}^{n_0/2}e_i\wedge\hat{e}_j\Bigl(\bigl(\nabla_{qe_i}\omega\bigr)(qe_j)+\frac{1}{2}\omega^2(qe_i,qe_j)\Bigr)\Bigr\}\,.
\end{align*}
Use Lemma \ref{Hodge} and $\Tr[\gamma_x^{\mathcal{F}}\omega^2(qe_i,qe_j)]=-d\theta(\gamma,\mathcal{F},\hF)_x(qe_i,qe_j)=0$.
\end{pf}
\begin{prop}\label{orient}
Let $M$ be oriented. For $n=\dim M$ even and $\gamma$
orientation preserving, or $n$ odd and $\gamma$
orientation reversing, one has for $x\in M$ and $t>0$:
\begin{equation*}
\Str\bigl[C_xk_t(\gamma,x)\bigr]=\Str^{\sigma}\bigl[C_xl_t(\gamma,x)\bigr]\,.
\end{equation*}
\end{prop} 
\begin{pf}
The proof of the corresponding non-equivariant result in \cite{BZ1} applies with the obvious modifications to the equivariant situation.
\end{pf}\\
\\
We have assumed up to this point that $M$ is oriented. Finally, for the evaluation of the asymptotic terms, we can remove this assumption since we can embed a $\gamma$-invariant neighbourhood of $x\in M^{\gamma}$ into a closed oriented Riemannian manifold such that the flat bundle $\F$ with its hermitian metric extends and so does the diffeomorphism $\gamma$, preserving the geometric data. The integrand $I(\gamma, \cdot)$, resp.~$I^{\sigma}(\gamma, \cdot)$, only depends on local geometric quantitites.\\
\\
{\it Proof} of {\bf Theorem \ref{Asymp}}. We consider a tubular neighbourhood $U \supset M^{\gamma}$. Let$\Psi_U$ be a cut-off function with support in $U$. Since the integral over $M\setminus U$ does not contribute asymptotically, we may write
\begin{align*}
\Str\bigl\{V\gamma\exp(-tH)\bigr\}&=\int_M\Str\bigl[Vk_t(\gamma,x)\bigr]\,|dx|\\
&\underset{t\rightarrow0}{\sim}\,\int_U\Str\bigr[Vk_t(\gamma,x)\bigr]\Psi_U(x)\,|dx|\\
&\underset{t\rightarrow0}{\sim}\,\int_{M^{\gamma}}\int_{\mathcal{N}_{x_0}}\Str\bigl[Vk_t(\gamma,\exp_{x_0}v)\bigr]\\
&\qquad\qquad\qquad\;\;\cdot\Psi_U(\exp_{x_0}v)j_U(v,x_0)\,|dv||dx_0|\,,
\end{align*}
where $j_U$ is the Jacobian of the exponential map, i.e.~$|dx|=j_U(x_0,v)|dv||dx_0|$, and $|dx_0|$ is the Riemannian density $M^{\gamma}$. We set $\phi(x_0,v)=\Psi(\exp_{x_0}v)j_U(v,x_0)$ and we get
\begin{align*}
\int_M\Str\bigl[k_t(\gamma,x)\bigr]|dx|
&\underset{t\rightarrow0}{\sim}\,\int_{M^{\gamma}} I(t,\gamma,\phi,x_0)\,|dx_0|\,,
\end{align*}
in particular 
\begin{align*}
\lim_{t\rightarrow0}\Str\bigl\{V\gamma\exp(-tH)\bigr\}=\int_{M^{\gamma}}I(\gamma,x_0)\,|dx_0|\,.
\end{align*}
In the same way we get 
\begin{align*}
\lim_{t\rightarrow0}\Str\bigl\{C\gamma\exp(-tH)\bigr\}=\int_{M^{\gamma}}I^{\sigma}(\gamma,x_0)\,|dx_0|\,.
\end{align*}
Assertion (2) of Theorem \ref{Asymp} now follows directly by substituting the formula for $I(\gamma,x_0)$ from Proposition \ref{I} and the definition of $e(TM^{\gamma},\nabla^{TM^{\gamma}})$:
\begin{align*}
\lim_{t\rightarrow0}\Str\bigl\{V\gamma\exp(-tH)\bigr\}=\int_{M^{\gamma}}\Tr[\gamma^{\mathcal{F}} V]e(TM^{\gamma},\nabla^{TM^{\gamma}})\,.
\end{align*}
Essentially the same calculation as in the non-equivariant case, cf.~\cite{BZ1}, yields
$$
\lim_{t\rightarrow0}\Str\bigl\{C\gamma\exp(-tH)\bigr\}=-\int\limits_{M^{\gamma}}\theta(\gamma,\mathcal{F},\hF)\wedge\tilde{e}^{\prime}(TM^{\gamma})\,,
$$
which proves Assertion (1).\hfill$\boxbox$

%%%%%%%%%%%%%%%%
% bibliography %
%%%%%%%%%%%%%%%%

\end{document}